\documentclass[titlepage,twoside,12pt]{article}
\usepackage{amssymb}
\usepackage{amsfonts}
\begin{document}

{\large \bf Reconsidering Conflict and Cooperation} \\
{\bf ~~~~~~~~~( revised version~ 1/7/5 )} \\ \\

{\it Elem\'{e}r E Rosinger \\ Department of Mathematics \\ University of Pretoria \\ Pretoria,
0002 South Africa \\
e-mail : eerosinger@hotmail.com} \\ \\

{\bf Abstract} \\

An analysis of several important aspects of competition or conflict in games, social choice
and decision theory is presented. Inherent difficulties and complexities in cooperation are
highlighted. These have over the years led to a certain marginalization of studies related to
cooperation. The significant richness of cooperation possibilities and the considerable gains
which my lie there hidden are indicated. Based on that, a reconsideration of cooperation is
suggested, as a more evolved form of rational behaviour. As one of the motivations it is shown
that the paradigmatic non-cooperative Nash equilibrium itself rests on a strong cooperation
assumption in the case of $n \geq 3$ players. \\ \\ \\

{\bf 0. Preliminary Remarks} \\ \\

Optimization has a long and important history in human endeavours. \\
In modern times, a significant moment occurred when the Newtonian laws of dynamics were found
to be expressible by the minimization of the integrals of corresponding Lagrangeans. This led,
as is well known, to the development of Variational Calculus. \\
Such minimization, or in general, optimization problems are, however, particular, and as such,
rather simple cases, since they involve a situation where one single conscious agent, a human
for instance, is facing what is usually called Nature. Indeed, two features are specific to
such a situation. Nature is supposed to have well set laws which act constantly and in the
same manner, regardless of the possible involvement of the conscious agent. Also, the
optimization pursued by the conscious agent has a clear cut and a priori well defined
criterion according to which is supposed to be accomplished. \\

Needless to say, there are other, more complex optimization situations in which several
conscious agents may become involved. Here we list {\it three} of such well known situations
which, so far, happen to exhaust a wide range of endeavours in optimization, see Luce \&
Raiffa. \\

{\it Games} are optimization situations in which two or more conscious and autonomous agents
are involved. A further feature in this case is that each such agent has a clear cut and a
priori well defined criterion according to which is supposed to act. \\

{\it Social, collective or group choice} is a situation in which one single conscious agent is
supposed to optimize the outcome for two or more beneficiaries. And again, each such
beneficiary is supposed to have a clear cut and a priori well defined criterion according to
which is supposed to be satisfied. \\

The third possibility is given by a {\it single decision maker with multiple and conflicting
objectives}. \\

At first sight, it may appear that games are the most difficult to deal with, since they
involve more than one conscious agent. Certainly, as we shall see, games are not easy to treat
in many cases of interest, and in fact, when approached in sufficient generality, they may
lead to algorithmically unsolvable problems. \\
However, as it turns out, social, collective or group choice is not a much simpler venture
either, even if there is only one conscious agent and the criteria which have to be satisfied
are clear cut and a priori well defined. Arrow's celebrated paradox, mentioned in the sequel,
gives a good measure of the difficulties involved in this regard. \\
Lastly, the situation of a single decision maker with multiple and conflicting objectives is
again not trivial, even if in this case he or she is the single beneficiary which has to be
satisfied. Indeed, the main difficulty in this case is in the fact that there is {\it no} -
and in general, there simply {\it cannot} be - a natural, canonical, universal way to
aggregate or synthesize a given set of multiple and conflicting objectives into one single
clear cut criterion. \\ \\

{\bf 1. Introduction} \\ \\

As mentioned, there are a number of well known ways interactions, and in particular,
competition or conflict, involving rational autonomous agents are modelled in mathematics,
among them by game theory, social, collective or group choice theory, and decision theory.
Obviously, such interactions, even when they involve competition or conflict, need not always
or necessarily be totally incompatible with certain forms of cooperation. And needless to say,
it is often that certain convenient outcomes cannot be obtained in any other ways, except by
cooperation. \\

Game theory, as initiated in its modern and systematic form by John von Neumann in the late
1920s, see von Neumann \& Morgenstern, centers around the {\it individual} rational agent and
aims to lead to a rational outcome when two or more such individuals interact in well defined
situations, and do so, however, without the interference of an overall arbitrating authority.
Social, collective or group choice theory also starts from individuals, yet its aim is to {\it
transcend} them to some extent by aggregating in suitable manners their given preferences, see
Luce \& Raiffa, Mirkin. \\
In a way at the other end of the spectrum involving rational autonomous agents, decision
theory aims, among others, to support in his or her decision one {\it single} decision maker
who faces a number of conflicting objectives, see Luce \& Raiffa, Bacharach \& Hurley,
Rosinger [1-4]. It follows that decision theory can be seen as variant of a one person game,
in which the player plays against Nature, the difficulty arising from the presence of several
conflicting objectives. \\

Here we can note that the very development of these mathematical theories and of the resulting
methods is in itself an act of rational behaviour, albeit on a certain {\it meta-level}, which
involves both the levels of the interest in developing the general concepts, axioms, theorems,
and so on, as well as the levels at which they are put to their effective uses in a variety of
relevant applications. \\

And yet, such a {\it meta-rational} behaviour appears to have mostly come to a halt when
dealing with {\it cooperation}. \\
Instead, the effort has rather been focused upon non-cooperative contexts, and rationality got
thus limited to them. A further aggravation of such a limitation upon rationality has come
from the fact that cooperation is so often being associated with, or even reduced to issues of
ethics, wisdom, philosophy, or on the contrary, to mere opportunistic coalitions based on
politics, or various other expediencies. \\

One of the few more prominent contributions to the study of cooperation has been the 1984 book
of Robert Axelrod, which however is limited to two person nonzero sum games, and it centers
around one of the simplest nontrivial such examples, namely, the celebrated game called the
Prisoner's Dilemma. This game was suggested in the early 1950s by Merrill Flood and Melvin
Drescher, and then formalized by Albert W Tucker, see section 3, or Axelrod, Rasmusen. \\

However, even such simplest nontrivial examples can show that cooperation itself is often but
a most natural matter of rationality, albeit manifested in forms which can often be far more
evolved than those encountered in non-cooperation. Indeed, in many non-cooperative situations
what can be obtained by the players involved proves to be significantly less than what may be
available through suitable cooperation. Thus the choice of cooperation need not at all be seen
as merely a matter which has to do with expediency, politics, wisdom, ethics, or morality.
Instead, cooperation proves to be one of the major {\it strategic assets} available in many
important situations. \\

And then the issue is simply the following : do we limit rationality and stop it before
considering cooperation in ways more adequate to its considerable depth and potential, or
instead, are we ready to try to be rational all the way ? \\

And the fact is that very few situations are of a nature in which competition or conflict is
total, thus, there cannot be a place for one or another form of cooperation. \\

In game theory, for instance, such a situation corresponds to the extreme and simplest case of
the two person zero sum game. \\
In rest, that is, in the vast majority of cases of games, the possibilities for cooperation
and its considerable possible rewards are always there. And then, all it takes is to {\it
extend} and {\it deepen} rationality, and thus find and develop suitable ways of cooperation.
And above all, to establish a context in which cooperation - which is a purely joint voluntary
matter - can indeed be relied upon. \\

In this regard it is instructive to see that, contrary to customary perception, even in the
celebrated paradigmatic {\it non-cooperative} Nash-equilibrium concept, and the corresponding
theorem, a {\it very strong cooperation} assumption is in fact essentially involved, when
there are n $\geq$ 3 players. And as it turns out, this cooperation assumption is indeed so
strong that it is simply {\it unrealistic} in practical situations. \\

In this way what we face is the alternative : \\

{\bf Non-cooperation} ~:~ Either we limit rationality and try to use it as much as possible
within the framework of non-cooperation, and only on occasion, and only as a second choice do
we consider and try cooperation as well. \\

Or \\

{\bf Cooperation} ~:~ Within an extended and deeper sense of rationality we create and
maintain a context in which due consideration given to cooperation can be relied upon. In
other words, we are {\it reliably cooperation minded}. \\

So far, as will be argued, when it comes to the development of the respective theories, the
first of the above two alternatives has mostly been taken, even if that may have happened
rather by default, that is, as a consequence of the extreme complexity of the phenomenon of
cooperation. \\
However, this limitation of rationality mostly to the study of non-cooperative games is not in
fact securing a significantly simpler or easier situation. Indeed, as follows from Binmore
[1-3], such games can still lead to complexities which are not algorithmically solvable. \\

In this way, having mostly avoided what appeared a theoretically difficult task in the second
above alternative, instead of it, the first alternative, which in fact is not less difficult,
was taken. Certainly, during the late 1940s and early 1950s when game theory knew a first
massive development, there was not much awareness about the possibility of the presence of the
type of deep difficulties which would more than three decades later be pointed out by Binmore.
On the contrary, during that first enthusiastic and major development period it was hoped that,
at last, game theory would offer mankind a theoretical model for dealing with all possible
conflicts among autonomous rational agents. And that hope was so deeply entrenched that Cold
War strategies at the time were suggested in the USA, based on game theoretical
considerations. \\

One way out, suggested in this study, and at this stage already with the benefit of knowing
about the message of Binmore, is to take the second above alternative, and do so in the
following manner. Now that it has been revealed that both the non-cooperative and cooperative
types of games can be of an extreme complexity from the point of view of theoretical
approaches, the focus can shift from the earlier attempts to construct comprehensive enough
theories for either of them, to finding a framework which is large enough to contain enough of
what is already known and it is important about games, and then, in such a framework, to
develop theories, methods, examples, applications, and so on, which have their own value and
interest, even if they may fall short of being ultimately or near ultimately comprehensive. \\

And clearly, once such a shift of focus is found appropriate, the {\it second alternative} is
the {\it natural} one. Indeed, this second alternative does obviously contain the first one as
a particular case. On the other hand, the first alternative has so far not been proven to
contain the second alternative. This is contrary to the hopes, or claims, of what once had
been the tempting "Nash Program", and which is mentioned in some more detail later. \\

Here it should be noted that cooperation is not limited to games alone. Indeed, in social,
collective or group choice, there may similarly exist a considerable scope for cooperative
interaction between the individual beneficiaries involved, an interaction which can, among
others, lead to a modification of some of their initial individual preferences, and also to a
readiness to accept an aggregations of such preferences which would otherwise be seen as being
less than fair to certain of the individual beneficiaries. Such modifications of individual
preferences can be arranged upon suitable side payments, similar to those in game theory. And
in view of the celebrated impossibility result of Kenneth J Arrow, discovered in 1950, see
section 4, such a cooperation may have the significant merit of being the one possible way to
avoid the necessity of a dictator. \\

At last, even in the case of one single decision maker, the fact that he or she faces all
alone his or her own conflicting objectives need not completely rule out the use of certain
cooperative type approaches. After all, cooperation can involve bargaining, and in the case of
one single decision maker with several conflicting objectives, he himself, or she herself may
end up as if bargaining with himself or herself. \\
In fact, certain forms of cooperation find a most appropriate context precisely within the
thinking of a single decision maker who exhibits a rational behaviour. \\

Of course, as mentioned, there are certain particular situations when cooperation is simply
not possible, like for instance, in the two person zero sum games, which in fact are the
simplest and most extreme instance of conflict in game theory.  On the other hand, in large
varieties of other interactions between autonomous rational agents there are vast and far from
sufficiently explored possibilities for cooperation. \\

As a rather relevant example in this respect it will be shown in section 2, that the
celebrated Nash-equilibrium concept and theorem itself, which is usually seen as a paradigm of
non-cooperation, must in fact assume, even if somewhat implicitly, a very strong form of {\it
cooperation}, in order to have any practical meaning and value at all in the case of $n \geq$
3 players. And then, paradoxically, the respective cooperation assumption is so strong that it
becomes unrealistic in practice. \\
Furthermore, the concept of Nash-equilibrium on which this result is based is highly {\it
unstable}, or {\it fragile}, when there are at least three players involved, players who are
not {\it cooperation minded}. \\
Another important criticism of both the Nash-equilibrium concept and of the corresponding
equilibrium result is that, as is well known, in case of cooperation, players can often obtain
significantly higher payoffs than those given by equilibrium strategies, see for instance the
Prisoner's Dilemma. \\

However, there is a rather weighty reason, explanation, and maybe also an excuse, for the fact
that meta-rational behaviour has so far mostly stopped before considering and dealing with
cooperation in ways proportionate to the significant advantages which are often available, and
so far, hidden there. Namely, as is well known even from common everyday experience,
cooperation will often involve considerable complications and difficulties. First of all, and
already on its most basic {\it conceptual} levels, cooperation proves to be an extremely
complex and rich phenomenon, which therefore cannot in any way be encompassed by a few general
definitions and mathematical models. This fact is, indeed, in sharp contrast with the
modelling of competition and conflict situations, where for instance in game theory, the so
called non-cooperative games, see (2.1) and (2.5) below, describe quite well, and in spite of
their manifest simplicity, a considerably large class of such situations. Second, one can only
talk about cooperation if one can {\it rely} on the respective agreements undertaken by the
autonomous agents involved. And the issue of such a reliance clearly depends on a variety of
complex factors which can easily be outside of the realms of convenient mathematical
modelling. \\

Yet, one should not forget that, just like in the case of competition and conflict, the
primary aspect of cooperation is {\it intent}, while the respective subsequent
conceptualizations, models and actions are only specific instances of manifestation and
expression of such an intent. Therefore, the primary issue is whether we intend to have
competition or conflict in a context which may preferably include cooperation as well, or on
the contrary, we intend, because of no matter what reasons, to relegate cooperation to a
secondary role, or even exclude it altogether. \\

And if we do not a priori intend to exclude cooperation, then we should be careful not to
allow that it is excluded merely by default, that is, due to the fact that it is in general
not so easy to deal with it, be it conceptually, or practically, and then, as a consequence,
it simply happens that we fail, avoid or decline to consider it. \\
And a readiness to pursue meta-rational behaviour beyond its present day limits confined
mostly to non-cooperation will then suggest that the intent to cooperate, and even more
importantly, {\it the intent to secure and keep up in the longer run a context suitable for
cooperation} is but a clearly rational behaviour, even if on a meta-level. \\

One of the most dramatic and clear cut theoretical manifestations of the major conceptual
difficulties involved related to cooperation is presented by the mentioned impossibility
result of Arrow. \\

The effect of the presence of such considerable difficulties related to cooperation has been
that the modelling of cooperation has not received enough attention, see Axelrod. An aim of
this study is to draw attention upon that fact, and also suggest certain ways to deal with
it. \\

Let us start with a few comments on three of the present day mathematical theories which deal
with interactions, among them competition, conflict and possibly cooperation, namely, see Luce
\& Raiffa, Rasmusen :

\begin{enumerate}

\item cooperative and non-cooperative games,

\item social, collective or group choice,

\item decision support systems for individual decision makers with conflicting objectives.

\end{enumerate}

Here we can mention that the way the customary division is made between cooperative and
non-cooperative games is in our view not quite appropriate, as we shall argue in section 2.
Indeed, it turns out that much of what is nowadays seen as a non-cooperative setup does in
fact rest, even if implicitly, on cooperative assumptions. For instance, as mentioned earlier,
and seen in section 2 in the sequel, the concept of Nash-equilibrium, and the respective Nash
theorem - both seen as paradigmatic for the non-cooperative games - turn out to be essentially
based on a certain very strong cooperative assumption in the case of $n \geq 3$ players. \\

Such and other similar prominent examples may provide some of the main points which one can
start dealing with related to a reconsideration of cooperation. And they may lead to certain
major practical consequences that may have been missed so far. \\
In particular, a good deal of what nowadays may be seen as non-cooperative games could
possibly be considered in the cooperative category, provided that suitable {\it extensions}
and {\it deepening} of the cooperation concept are employed. \\

Needless to say, such a view runs against, among others, the mentioned Nash Program of the
early 1950s, see Nash [1,3], which tried to achieve the opposite, by reducing cooperative
games to non-cooperative ones. \\

Certainly, such a program - or for that matter, one aiming the other way round - comes quite
likely from deeper and more general views about the nature of possible interaction among
autonomous rational agents, than the views customarily encountered. \\

As is well known, Nash himself tended to see such interactions as being mainly moved by
competition, conflict, and so on, rather than by cooperation, see Nasar. \\
On the other hand, at that time, the older and much more experienced John von Neumann, who was
in fact the originator of modern game theory, considered that there was a major and urgent
need in economics and other important human ventures involving strategic thinking for the
introduction of rational approaches to the respective variety of human interactions involved.
And clearly, the very attempts to rationalize approaches to competition, conflict, and so on,
rather tend to mollify, than prioritize them. Thus cooperation does not become something to be
avoided, but rather to be considered, and made use of, whenever possible. Consequently, von
Neumann chose to devote a large and important role to various forms of cooperation which, from
theoretical point of view, could be modelled and dealt with at the time, see Neumann \&
Morgenstern. \\

Rather independently during the same period came the famous 1950 paper of Arrow on the
impossibility of setting up in general a social choice function in the absence of a dictator.
As we shall see, however, this and the subsequent related developments in what is nowadays
called social, collective or group choice theory, can naturally be seen as stressing too much
the opposition between cooperation and non-cooperation, thus leading to a particular and often
less than welcome solution method, namely, the use of a dictator. On the other hand, by
extending and deepening the concept of cooperation one can attain suitable aggregations of
social, collective or group choice, and do so without the need for dictators. \\

Let us now mention in short a history of these and related events. A more detailed history of
game theory can be found in Luce \& Raiffa, or Walker, while for a view of the background to
group choice and decisions one can consult Bacharach \& Hurley, respectively, Mirkin. \\

The first major result in game theory was obtained by John von Neumann in his 1928 paper. This
is the famous Min-Max theorem about two player zero sum games. Such games involve the smallest
possible nontrivial number of players, and the most extreme possible conflict among them, in
which what one player wins, the other one must lose, thus the sum of what is won and lost is
always zero. Clearly, in such a game there is no any way available for cooperation. \\

Here it is important to note the following. During the period around 1928, when von Neumann
was only 25 years old, he was involved in at least two other major ventures, namely, the
foundation of set theory, and the foundation of quantum mechanics. And in both of them he made
most important and lasting contributions. In this way, von Neumann's involvement in game
theory during that period can be seen as reflecting the special, if not in fact, fundamental
importance he happened to attribute to it. And indeed, he saw it as being the first ever
systematic and rigorous theoretical approach to a {\it rational} management of conflict or
competition between two or more conscious agents. \\

Problems of optimization had been considered earlier as well. After all, much of the practical
human endeavours have always been pursuing one or another from of optimization. \\
One major example, starting with the 1800s, had been given by variational calculus which
proved to lead to an equivalent formulation of Newtonian mechanics. \\
However, such problems could be seen as a game with one single conscious and rational player
who was playing against Nature. \\

But now the task facing von Neumann was to be able to build an appropriate theory for
conflicts and competitions between two or more conscious agents, assuming that they were
firmly and reliably grounded in rationality. \\
It should be remembered in this regard that von Neumann happened to grow up in the Empire of
Austria-Hungary, and did so during the disastrous years of World War I, and its immediate
aftermath. And just like the well known philosopher Karl Popper, of the same generation and
social background, von Neumann was much influenced by views dominant in the post World War I
years. Views according to which World War I - called at the time The Great War -  was seen as
nothing else but a systematic and catastrophic, even if rather trivial succession of failures
of rationality on the part of the elites running the Western powers. \\

The extent to which von Neumann gave a special priority to the development of game theory is
further illustrated by his activity during the next one and a half decade, till the publishing
in 1944 of his joint book with Morgenstern, entitled "Theory of Games and Economic Behaviour".
Indeed, during the years of World War II, von Neumann was heavily involved in supporting the
American war effort and doing so in a variety of ways. Consequently, at the time, he did very
little theoretical research. And yet, he considered it important enough to dedicate time to
game theory, and complete the mentioned book of over 600 pages, which is the first ever
systematic and detailed presentation of that theory. It should also be mentioned that the
theory in that book is due solely to von Neumann, and most of it, except for his Min-Max
theorem of 1928, was developed by him in the period leading to its first publication in 1944.
Morgenstern was an economist, and his contribution to the book consisted in the connections
between game theory and economic behaviour. In this way, that book can in fact be seen as a
massive research monograph - and in fact, the first one - in game theory. \\

The importance attributed to game theory continued after World War II as well. And it was due
to a good extent to the interest manifested in it at the RAND Corporation, a most influential
California based think tank at the time, which was heavily involved, among others, in
strategic studies related to the emerging Cold War. \\

As it happens, Emile Borel initiated in the early 1920s the study of certain well known card
games which were related to the two person zero sum games. However, he did not obtain the
respective major result, namely, the Min-Max theorem, and in fact, he assumed that such a
theorem was in general false. Later, in 1934, the well known statistician R A Fisher was also
involved in a study of two person zero sum games, without however obtaining the major Min-Max
theorem, see Luce \& Raiffa. \\

Then starting in 1950, John Nash, who at the time was 22 years old, published his fundamental
papers, Nash [1,3], on equilibrium in n-person non-cooperative games, and his main result was
a significant and two fold extension of the Min-Max theorem. Indeed, unlike von Neumann's
Min-Max, Nash managed to show the existence of equilibrium not only in the case of two players
and a zero sum game, but for an arbitrary finite number of players, and for games with an
arbitrary sum. \\
What appeared to be similar with the earlier Min-Max theorem, was the non-cooperative nature
of the result of Nash. That similarity, however, will be shown to be but illusory, see section
2. \\
Indeed, von Neumann's Min-Max theorem being about two players and with zero sum, it does not
have any other chance but to be non-cooperative. \\
On the other hand, both the concept of Nash-equilibrium and the respective Nash theorem are,
as mentioned, essentially based on a very strong cooperative type assumption in the case of
$n \geq$ 3 players. \\

By the way, here it should be mentioned that Nash also published important result in
cooperative games, Nash [2,4], see section 3. \\

However, as it happened, his result, which massively extended the Min-Max theorem to arbitrary
number of players and nonzero sum games, and was later in the 1990 to earn him a Nobel Prize
in Economics, has ever since been seen as essentially belonging to non-cooperative games. \\
To a certain extent, such an interpretation is not so surprising due to the following two
facts. First, the Nash result on the existence of an equilibrium in mixed strategies is an
obvious extension of the Min-Max theorem of von Neumann, and the latter, as mentioned, is
indeed about games which are outside of any possible cooperation. Second, as long as one is
limited to the usual, and thus narrow concepts of cooperation, the result of Nash will be seen
as falling outside of the cooperative framework. \\

As we shall show in section 2, however, such an interpretation can only hold if the usual, and
indeed narrow concepts of cooperation are considered. On the other hand, as the very concept
of equilibrium in the Nash result implies it, that result can have any practical meaning and
value at all, and do so beyond its particular two person zero sum Min-Max case, only if the
respective $n \geq$ 3 players do accept - even if implicitly - certain very strong {\it
additional common rules} of behaviour. Thus in the case of $n \geq 3$ players, they {\it must}
end up by cooperating very strongly, even if in ways other, and more deep, than those
according to the usual views of cooperation. \\
Indeed, the kind of cooperation needed in order to enable the Nash equilibrium concept and
result to function at all proves to be particularly {\it strong}. So strong, in fact, as to
render it {\it unrealistic} in practice. \\

Due to the reputation of von Neumann, the interest showed by him in game theory led in the
late 1940s and early 1950s to a considerable status for that theory among young mathematicians
at Princeton, see Nasar. That status was further enhanced by results such as those of Nash and
a number of other mathematicians. \\
There was also at the time a significant interest in game theory outside of academe. As
mentioned for instance, the influential RAND Corporation was conducting studies in political
and military strategy which were modelled mathematically by a variety of games. \\

As it happened, however, soon after, certain major setbacks were experienced. First, and
within game theory itself, was the fact that in the case of n-person games, even for $n \geq$
3 moderately large, there appeared to be serious {\it conceptual} difficulties related to
reasonable definitions of solution. Indeed, too many such games proved not to have solutions
in the sense of a variety of solution concepts, concepts which each seemed to be natural, see
Luce \& Raiffa, Owen, Vorob'ev, Rasmusen. \\
Later, the nature and depth of these conceptual difficulties got significantly clarified. For
instance, in Binmore [1-3] it was shown that there are no Turing machines which could compute
general enough games. In other words, solving games is not an algorithmically feasible
problem. \\

The second major trouble came from outside of game theory, and it is not quite clear whether
at the time it was soon enough appreciated by game theorists themselves, with respect to its
possible implications about the fundamental difficulties in formulating appropriate concepts
of solutions in games. Namely, as mentioned, Arrow showed that a set of individual preferences
cannot in general, and under reasonable conditions, be aggregated into one joint preference,
unless there is a dictator who can impose such a joint preference. \\

Strangely enough, this result of Arrow about difficulties in the aggregation of a number of
individual preferences was not completely new or unknown. Indeed, it was in fact extending and
deepening the earlier known, so called, {\it Voter's Paradox}, mentioned by the Marquis de
Condorcet, back in 1785, see section 4, or Mirkin. \\

In subsequent years, following the 1950s, developments in game theory lost much of their
momentum. In more recent years, applications of game theory gained the interest of economists,
and led to a number of new developments in economic theory. An indication of such developments
was the founding in 1989 of the journal Games and Economic Behavior. \\

At the same time certain studies of competition, conflict, and so on, were taken up by the
developments following Arrow's fundamental paper, and leading to social, collective or group
choice theory, among others. In this regard, a considerable literature has been developed. \\
Decision theory got also involved in such studies involving certain specific instances of
competition, conflict or cooperation, related to problems of optimization, see Rosinger
[1-4]. \\

In games, or in social, collective or group choice one has many autonomous players,
participants or agents involved, each of them with one single objective, namely, to maximize
his or her advantage which usually is defined by a scalar, real valued utility function. And
by the early 1950s it became clear enough that such a situation would not be easy to handle
rationally, even on a conceptual level. \\
On the other hand, a main objective of decision theory is to enable one single decision maker
who happens to have several different, and usually, quite strongly conflicting objectives. And
in view of Arrow's result, such a situation may appear to be more easy to deal with since the
single decision maker can anyhow function as a "dictator" for himself or herself. \\
Yet as seen in section 6, being even such a special case of "dictator" will present its rather
difficult problems. \\ \\

{\bf 2. The Nash Equilibrium and Theorem} \\ \\

The usual way an n-person non-cooperative game in terms of the players' {\it pure strategies}~
is defined is by \\

(2.1) \quad $ G ~=~ ( P,~ ( S_i ~|~ i \in P ),~ ( H_i ~|~ i \in P ) ) $ \\

Here $P$ is the set of $n \geq 2$ players, and for every player $i \in P$, the finite set
$S_i$ is the set of his or her pure strategies, while $H_i : S \longrightarrow \mathbb{R}$ is
the {\it payoff} of that player. Here we denoted by \\

(2.2) \quad $ S ~=~ \prod_{i \in P}~ S_i $ \\

the set of all possible {\it aggregate pure strategies} ~$s = ( s_i ~|~ i \in P) \in S$
generated by the independent and simultaneous individual strategy choices $s_i$ of the players
$i \in P$. \\

The game proceeds as follows. Each player $i \in P$ can freely choose an individual strategy
$s_i \in S_i$, thus leading to an aggregate strategy $s = ( s_i ~|~ i \in P) \in S$. At that
point, each player $i \in P$ receives the payoff $H_i ( s )$, and the game is ended. We assume
that each player tries to maximize his or her payoff. \\ \\

{\bf Remark 1} \\

The usual reason the games in (2.1) are seen as non-cooperative is as follows. Each of the
$n \geq 2$ players $i \in P$ can completely independently of any other player in $P$ choose
any of his or her available strategies $s_i \in S_i$. And the only interaction with other
players happens on the level of payoffs, since the payoff function $H_i$ of the player $i \in
P$ is defined on the set $S$ of aggregated strategies, thus it can depend on the strategy
choices of the other players. \\
However, as we shall see in Remarks 2 - 4 below, in the case of $n \geq 3$ players, this
independence of the players is only apparent, when seen in the framework the concept of
Nash-equilibrium, and the corresponding celebrated Nash theorem.

\hfill $\Box$ \\ \\

Before considering certain concepts of equilibrium, it is useful to introduce some notation.
Given an aggregate strategy $s = ( s_i ~|~ i \in P ) \in S$ and a player $j \in P$, we denote
by $s_{-j}$ what remains from $s$ when we delete $s_j$. In other words $s_{-j} = ( s_i ~|~
i \in P \setminus \{ j \} )$. Given now any $s^{\,\prime}_j \in S_j$, we denote by $( s_{-j},
s^{\,\prime}_j )$ the aggregate strategy $( t _i ~|~ i \in P ) \in S$, where $t_i = s_i$, for
$i \in P \setminus \{ j \}$, and $t_j = s^{\,\prime}_j$, for $i =j$.  \\

For every given player $j \in P$, an obvious concept of {\it best strategy} $s^*_j \in S_j$ is
one which has the {\it equilibrium} property that \\

(2.3) \quad $ H_j ( s_{-j}, s^*_j ) ~\geq~ H_j ( s_{-j}, s_j ),~~~ \mbox{for all}~~
                                                                s \in S,~ s_j \in S_j $ \\

Indeed, it is obvious that any given player $j \in P$ becomes completely independent of all
the other players, if he or she chooses such a best strategy. However, as it turns out, and is
well known, Rasmusen, very few games of interest have such strategies. \\
Consequently, each of the players is in general {\it vulnerable} to the other players, and
therefore must try to figure out the consequences of all the possible actions of all the other
players. \\

Furthermore, even when such strategies exist, it can easily happen that they lead to payoffs
which are {\it significantly lower} than those that may be obtained by suitable cooperation. A
good example in this regard is given by game called the Prisoner's Dilemma, see section
3. \\ \\

{\bf Remark 2} \\

It is precisely due to the mentioned vulnerability of players, which is typically present in
most of the games in (2.1), that there may arise an interest in {\it cooperation} between the
players. A further argument for cooperation comes from the larger payoff individual players
may consequently obtain. A formulation of such a cooperation, however, must then come in {\it
addition} to the simple and general structure present in (2.1), since it is obviously {\it not}
already contained explicitly in that structure.

\hfill $\Box$ \\ \\

Being obliged to give up in practice on the concept of best strategy in (2.3), Nash suggested
the following alternative concept which obviously is much {\it weaker}. \\ \\

{\bf Definition ( Nash )} \\

An aggregate strategy $s^* = ( s^*_i ~|~ i \in P) \in S$ is called a {\it Nash-equilibrium},
if no single player $j \in P$ has the incentive to change {\it all alone} his or her strategy
$s^*_j \in S_j$, in other words, if

\bigskip
(2.4) \quad $ H_j ( s^* ) ~\geq~ H_j ( s^*_{-j}, s_j ),~~~ \mbox{for all}~~ j \in P,~ s_j \in S_j $ \\ \\

{\bf Remark 3} \\

Clearly, the Nash-equilibrium only considers the situation when {\it never more than one
single} player does at any given time deviate from his or her respective strategy. Therefore,
the Nash-equilibrium concept is {\it not} able to deal with the situation when there are $n
\geq$ 3 players, and at some moment, more than one of them deviates from his or her
Nash-equilibrium strategy. \\
Needless to say, this fact renders the concept of Nash-equilibrium {\it unrealistically
particular}, and as such, also {\it unstable} or {\it fragile}. \\

Furthermore, that assumption has a manifestly, even if somewhat implicitly and subtly, {\it
cooperative} nature. \\

Above all, however, the larger the number $n \geq$ 3 of players, the {\it less realistic} is
that assumption in practical cases. \\

It is obvious, on the other hand,  that when there are $n \geq 3$ players, in case at least
two players change their Nash-equilibrium strategies, the game may open up to a large variety
of other possibilities in which some of the players may happen to increase their payoffs. \\

Therefore, when constrained within the context of the Nash-equilibrium concept, the game
becomes {\it cooperative} by {\it necessity}, since the following {\it dichotomy} opens up
inevitably : \\

\newpage

Either

\begin{itemize}

\item (C1)~~~ All the players agree that never more than one single \\
          \hspace*{1.3cm} player may change his or her Nash-equilibrium strategy,

\end{itemize}

Or

\begin{itemize}

\item (C2)~~~ Two or more players can set up one or more coalitions, \\
          \hspace*{1.3cm} and some of them may change their Nash-equilibrium \\
          \hspace*{1.3cm} strategies in order to increase their payoffs.

\end{itemize}

Consequently, what is usually seen as the essentially non-cooperative nature of the game (2.1),
turns out, when seen within the framework of the Nash-equilibrium concept, to be based - even
if tacitly and implicitly - on the {\it very strong cooperative} assumption (C1) in the above
dichotomy. \\
On the other hand, in case (C1) is rejected, then the game falls out of the Nash-equilibrium
framework, and thus it opens up to the wealth of possibilities under (C2), which among others,
can contain a large variety of possible ways of cooperation. \\

In this way, both the Nash-equilibrium concept and the Nash theorem on the existence of the
respective equilibrium in mixed strategies are highly {\it unstable} or {\it fragile} when
there are 3 or more players involved. \\

Also, similar with the best strategies in (2.3), with the Nash-equilibrium
strategies as well it can happen that they lead to payoffs which are significantly lower than
those that may be obtained by suitable cooperation.

\hfill $\Box$ \\ \\

As in the particular case of (2.1) which gives the von Neumann Min-Max theorem on two person
zero sum games, so with the weakened concept of Nash-equilibrium in (2.4), such an equilibrium
will in general not exist, unless one embeds the {\it pure strategy} game (2.1) into its
extension given by the following {\it mixed strategy} game \\

(2.5) \quad $ \mu G ~=~ ( P,~ ( \mu S_i ~|~ i \in P ),~ ( \mu H_i ~|~ i \in P ) ) $ \\

Here, for $i \in P$, the set $\mu S_i$ has as elements all the probability distributions
$\sigma_i : S_i \longrightarrow [0, 1]$, thus with $\Sigma_{s_i \in S_i}~ \sigma_i ( s_i ) =
1$. Let us now denote \\

(2.6) \quad $ \mu S ~=~ \prod_{i \in P}~ \mu S_i $ \\

Then for $i \in P$ we have the payoff function $\mu H_i : \mu S \longrightarrow \mathbb{R}$
given by \\

(2.7) \quad $ \mu H_i ( \sigma ) ~=~  \Sigma_{s \in S}~ \sigma ( s ) H_i ( s ) $ \\

where for $\sigma = ( \sigma_i ~|~ i \in P ) \in \mu S$ and $s = ( s_i ~|~ i \in P ) \in S$,
we define $\sigma ( s ) = \prod_{i \in P}~ \sigma_i ( s_i )$. \\

Now the definition (2.4) of Nash-equilibrium for pure strategy games (2.1) extends in an
obvious manner to the mixed strategy games (2.5), and then, with the above we have, see
Vorob'ev  \\ \\

{\bf Theorem ( Nash) } \\

The mixed strategy extension $\mu G ~=~ ( P,~ ( \mu S_i ~|~ i \in P ),~ ( \mu H_i ~|~ i
\in P ) )$ of every pure strategy game $G ~=~ ( P,~ ( S_i ~|~ i \in P ),~ ( H_i ~|~ i
\in P ) )$, has at least one Nash-equilibrium strategy. \\ \\

{\bf Remark 4} \\

Obviously, what was mentioned in Remarks 2 and 3 related to the inevitability of cooperation when the {\it pure
strategy} games (2.1) are considered within the framework Nash equilibrium, will also hold for the {\it mixed strategy}
games (2.5), and thus as well for the above theorem of Nash. \\ \\

{\bf 3. Usual Cooperation Models} \\ \\

The idea behind the Nash Program to reduce cooperative games to non-cooperative ones seems at
first quite natural. Indeed, in its very essence, a game means that, no matter what the rules
of the game are, each player must nevertheless remain with a certain {\it residual freedom} to
act within those rules, and be able to do so {\it independently} of the other players involved.
Therefore, it may appear that if we only concentrate on that freedom and independence, then
within that context one can see the game as non-cooperative, that being one of the usual ways
to understand the very meaning of freedom and independence. \\
Indeed, even if one cooperates, one is still supposed to be left in a game with a certain
freedom and independence. Thus it may still appear that, after subtracting all what is due to
the rules of the game and to one's possible cooperation, one is still supposed to remain with
a certain freedom and independence. \\

Put in a simple formula, the Nash Program may appear as the statement

\begin{quote}

game ~~$\equiv$~~ non-cooperative game~ ( modulo~ cooperation )

\end{quote}

According to Nash himself, it could be possible to express all communication and bargaining in
a cooperative game in a formal manner, thus turn the resulting freedom and independence of the
players into moves in an {\it extended} non-cooperative game, in which the payoffs are also
extended accordingly. Since such a program has never been fully implemented in all its details
and only its ideas were presented, its criticism must unavoidably remain on the same level,
namely, of ideas. However, a certain relevant and well tested objection can be made
nevertheless, see McKinsey [p. 359] :

\begin{quote}

" It is extremely difficult in practice to introduce into the cooperative games the moves
corresponding to negotiations in a way which will reflect all the infinite variety permissible
in the cooperative game, and to do this without giving one player an artificial advantage
( because of his having the first chance to make an offer, let us say )."

\end{quote}

What is lost, however, in such a view as the Nash Program is that an {\it appropriate
voluntary and mutual limitation} of one's freedom and independence, in order to implement a
{\it cooperation} can {\it significantly change} the payoffs, and thus it can offer to players
an increase in their payoffs, an increase which simply {\it cannot} be attained in any other
non-cooperative way. And this is after all the point in {\it cooperation}. This is the basic
fact which makes cooperation useful, and thus interesting all in itself, and in particular,
irreducible. \\

On the other hand, precisely to the extent that the above objection in McKinsey is valid
related to the Nash Program, and all subsequent experience points to its validity, the very
same objection touches essentially as well on any attempt to reconsider cooperation, and do so
in more formal ways. \\

And then, the main issue when reconsidering cooperation is not so much a renewed effort in
trying to formalize its myriads of possibilities, but rather, in first deciding to create and
maintain a reliable context of cooperation, and then starting from there, to develop in
various more typical and important instances certain appropriate theories as well. \\

We shall now return to illustrate in a few simple examples the mentioned advantages of
cooperation. And in two person nonzero sum games cooperation is easier to consider since there
is only one way for it, namely, by the involvement of both, thus of all, players. In other
words, in terms of coalitions, see von Neumann \& Morgenstern, in a two person game there is
only one possible coalition. Therefore we start here with this simple case by recalling in its
main features the way such a cooperation was addressed in Nash [2], see also Luce \& Raiffa. A
more general approach to the issue of cooperation is presented in section 5. \\

One of the most simple and dramatically clear examples which show that there can be a most
important point in cooperation even in the case of such a simple situation as a two person
nonzero sum game is given by the Prisoner's Dilemma which we present now in short. \\

Assume therefore that with the notation in (2.1), (2.2), we have $P ~=~ \{~ 1, 2 ~\}$ and $S_1
~=~  S_2 ~=~ \{~ 1, 2 ~\}$, then $H_1$ and $H_2$~can be given by the double matrix

$$ \begin{array}{l}
                               H_1(~1,~1),~ H_2(~1,~1) ~~~~~~ H_1(~1,~2),~ H_2(~1,~2) \\ \\

                               H_1(~2,~1),~ H_2(~2,~1) ~~~~~~ H_1(~2,~2),~ H_2(~2,~2)
      \end{array} $$

\medskip
In the case of the Prisoner's Dilemma this double matrix is given by \\

$$ \begin{array}{l}
                               ~~~~~~~~~~~~~~~~~~~~~~deny ~~~~~~~~~~~~~~~~~~~~~confess \\ \\

                               deny ~~~~~~~~~~~~- 1,~ - 1 ~~~~~~~~~~~~~~~~~~- 20,~~~~ 0 \\ \\

                               confess ~~~~~~~~~~~~ 0,~ - 20 ~~~~~~~~~~~~~~~~~- 8,~ - 8
     \end{array} $$ \\

\medskip
and has the following meaning. The two prisoners are not allowed to communicate with one
another. If both confess to the crime they are jointly accused of, then both get 8 years in
prison. If both deny the crime, then both get only 1 year in prison. Finally, if one denies
and the other confesses, then the one denying it gets 20 years in prison, while the one
confessing it is set free. \\
Clearly, in this case we have the respective strategy sets \\

$~~~~~~ S_1 ~=~ S_2 ~=~ \{~ 1 = deny,~ 2 = confess ~\} $ \\

and the game is symmetric with respect to the two players, as far as their individual
strategies and payoffs are concerned. \\

It is now obvious that for each player the best strategy in the sense of (2.3) is to confess
to the crime. Indeed, in such a case one will at worst get 8 years in prison, and avoids a 20
years sentence which could fall upon the one who denies the crime, while the other confesses
to it. However, precisely because of the symmetry of the situation, if both choose their best
strategies and confess, then both will stay 8 years in prison, while had they cooperated and
thus both denied the crime, they would have each escaped with only 1 year in prison. \\

This game thus shows immediately the following facts : \\

1. The best strategies in the sense of (2.3) do not always lead to the best payoffs. \\

2. Cooperation can give better payoffs than the best strategies. \\

3. The best payoffs are not always available, not matter whether one cooperates or not. \\

Of course, the Prisoner's Dilemma can model a large variety of other two person nonzero sum
symmetric games which may have other payoff functions. Also, $n \geq 3$ person variants of the
Prisoner's Dilemma have been studied, see Axelrod. \\

Needless to say, with respect to cooperation a large variety of complex issues arise even in
such a simple case as that of the Prisoner's Dilemma, see Axelrod. To mention just one such
issue, let us note the important difference it makes whether the game is played only once, a
fixed and known finite number of times, or a number of times which is not known from the
beginning. And clearly, in the last case, for instance, the relevance of cooperation can only
increase significantly. \\

Let us no turn to the {\it axiomatic} approach to bargaining, and thus to one form of
cooperation suggested by Nash. This approach is quite typical for the period of late 1940s and
early 1950s in game theory, when the basic concepts of best strategies, respectively,
equilibrium were seen as being strongly related to certain principles of fairness, plus
possibly some suitable mathematical conditions, which together would lead to the proper, if
not even unique, definition of the respective concepts, see Rasmusen. \\

The axioms proposed by Nash were meant to be so natural, fair and reasonable, as to impel the
two players to accept them, thus establishing the cooperation between them. In this way, such
axioms belonged to the {\it meta-rational} level of approach to game theory, and as such, they
were a good example for extending and deepening the concept of cooperation. As it turned out,
however, no matter how fair, reasonable and natural the axioms of Nash were, very simple
examples proved that in certain applications they would lead to strange, questionable, if not
in fact, unacceptable outcomes. In this way, once again it became obvious - even if at the
time Nash suggested them, not as clearly as from the later study of Binmore - that game theory,
and above all, the issue of cooperation, are extremely complex phenomena. \\

The versions of the two person nonzero sum games which Nash considered for cooperation or
bargaining are given by a {\it convex, closed} and {\it bounded} set $H$ in $\mathbb{R}^2$
which describes all the possible outcomes for the pairs of payoffs of the two players $~i \in
P ~=~ \{~ 1, 2 ~\}$, both of whom are supposed to want to maximize their respective payoffs.

\newpage

\bigskip
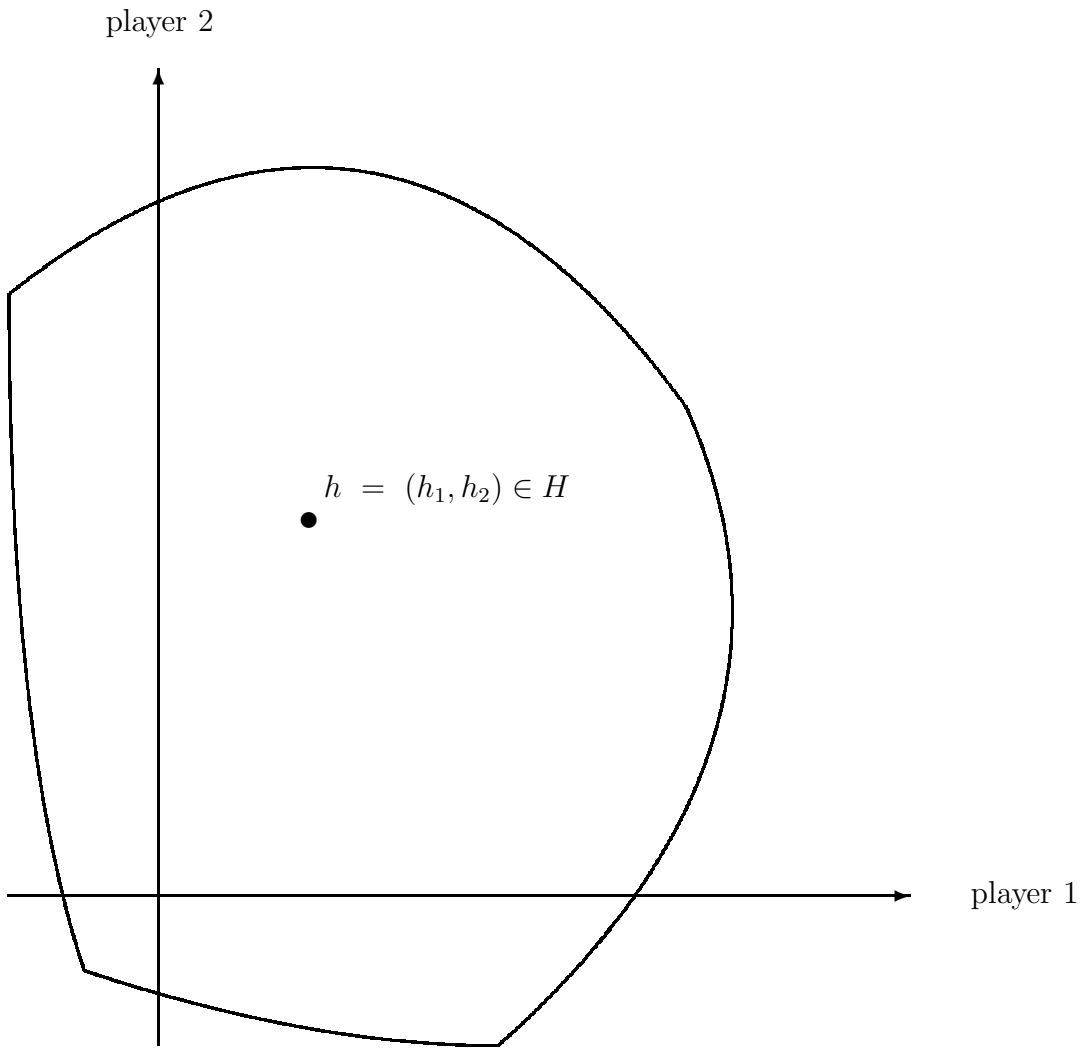
\begin{math}
\setlength{\unitlength}{1cm}
\thicklines
\begin{picture}(15,14)
\put(2,0){\vector(0,1){13}}
\put(1.3,13.5){$\mbox{player~2}$}
\put(0,2){\vector(1,0){12}}
\put(12.8,1.9){$\mbox{player~1}$}
\qbezier(0,10)(5,14)(9,8.5)
\qbezier(9,8.5)(11,4)(6.5,0)
\qbezier(6.5,0)(4,0)(1,1)
\qbezier(1,1)(0,4)(0,10)
\put(4,7){\circle*{0.2}}
\put(4.2,7.3){$h ~=~ ( h_1, h_2 ) \in H$}
\end{picture}
\end{math} \\

\medskip
Namely, if the outcome of the game is the point $~h ~=~ ( h_1, h_2 ) \in H$ above, then player
1 receives a payoff of $h_1$, while player 2 receives a payoff of $h_2$.\\ \\

{\bf Remark 5} \\

Clearly, since the set $H$ of all the possible payoffs is in general an infinite set of points,
thus of possible outcomes, the corresponding two person games cannot be in the pure strategy
form (2.1), but rather in the mixed strategy from (2.5). However, it is easy to see that {\it
not} all non-cooperative games in mixed strategy form (2.5) lead to a convex set $H$ of
possible payoffs. Indeed, let us consider the following two person nonzero sum non-cooperative
game in pure strategy form given by the double matrix \\

$$ \begin{array}{l}

                               ~~~~~~~~~~~~~~ 2,~~~~~~ 1 ~~~~~~~~~~~~~~ - 1,~~ - 1 \\ \\

                               ~~~~~~~~~~~ - 1,~~~ - 1 ~~~~~~~~~~~~~~~~~ 1,~~~~~ 2
     \end{array} $$

\medskip
Then the corresponding mixed strategy form of this game will have the set $H$ of possible
payoffs is given by, see Luce \& Raiffa [p. 93]

\newpage

\bigskip
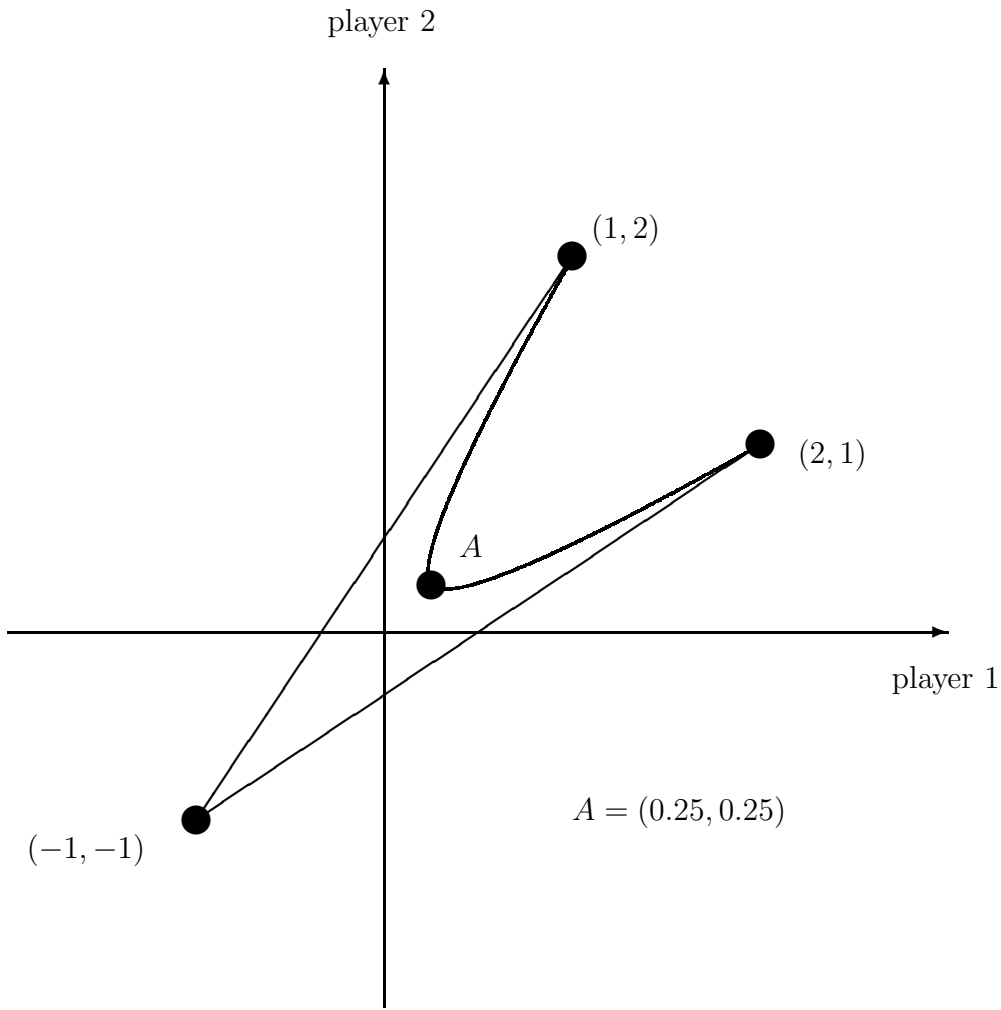
\begin{math}
\setlength{\unitlength}{2.5cm}
\thicklines
\begin{picture}(5,5.5)
\put(2,0){\vector(0,1){5}}
\put(1.7,5.2){$\mbox{player~2}$}
\put(0,2){\vector(1,0){5}}
\put(4.7,1.7){$\mbox{player~1}$}
\put(0.1,0.8){$( - 1, - 1 )$}
\put(1,1){\circle*{0.15}}
\put(1,1){\line(2,3){2}}
\put(3.1,4.1){$( 1, 2 )$}
\put(3,4){\circle*{0.15}}
\put(1,1){\line(3,2){3}}
\put(4.2,2.9){$( 2, 1 )$}
\put(4,3){\circle*{0.15}}
\qbezier(3,4)(2.1,2.4)(2.25,2.25)
\qbezier(2.25,2.25)(2.4,2.1)(4,3)
\put(2.25,2.25){\circle*{0.15}}
\put(2.4,2.4){$A$}
\put(3,1){$A =( 0.25, 0.25 )$}
\end{picture}
\end{math} \\

\medskip
which is obviously not convex. Therefore, there are plenty of two person nonzero sum
non-cooperative games which, when taken even in their mixed strategy form, cannot be contained
in the Nash bargaining or cooperative model. \\

However, even here, the advantages of cooperation already show up in an obvious manner. For
instance, in the above non-cooperative game, the payoff pair $( 1.5, 1.5 )$ is {\it not}
available to the two players even in the mixed strategy form, since it does not belong to the
nonconvex set $H$. \\
Yet, by cooperating, the two players can easily reach this point $( 1.5, 1.5 )$ which will
give each of them the payoff 1.5. Indeed, this can be done if they agree to go only for the
payoffs $( 2, 1 )$ and $( 1, 2 )$, and do so with equal frequency. \\

In general, it is easy to show that by cooperation, in every two person mixed strategy game
(2.5) it is possible to {\it extend} the initial and possibly nonconvex set $H$ of payoff
outcomes to its {\it convex closure} $H^{\#}$, which in the case of the above game will look
as follows

\bigskip
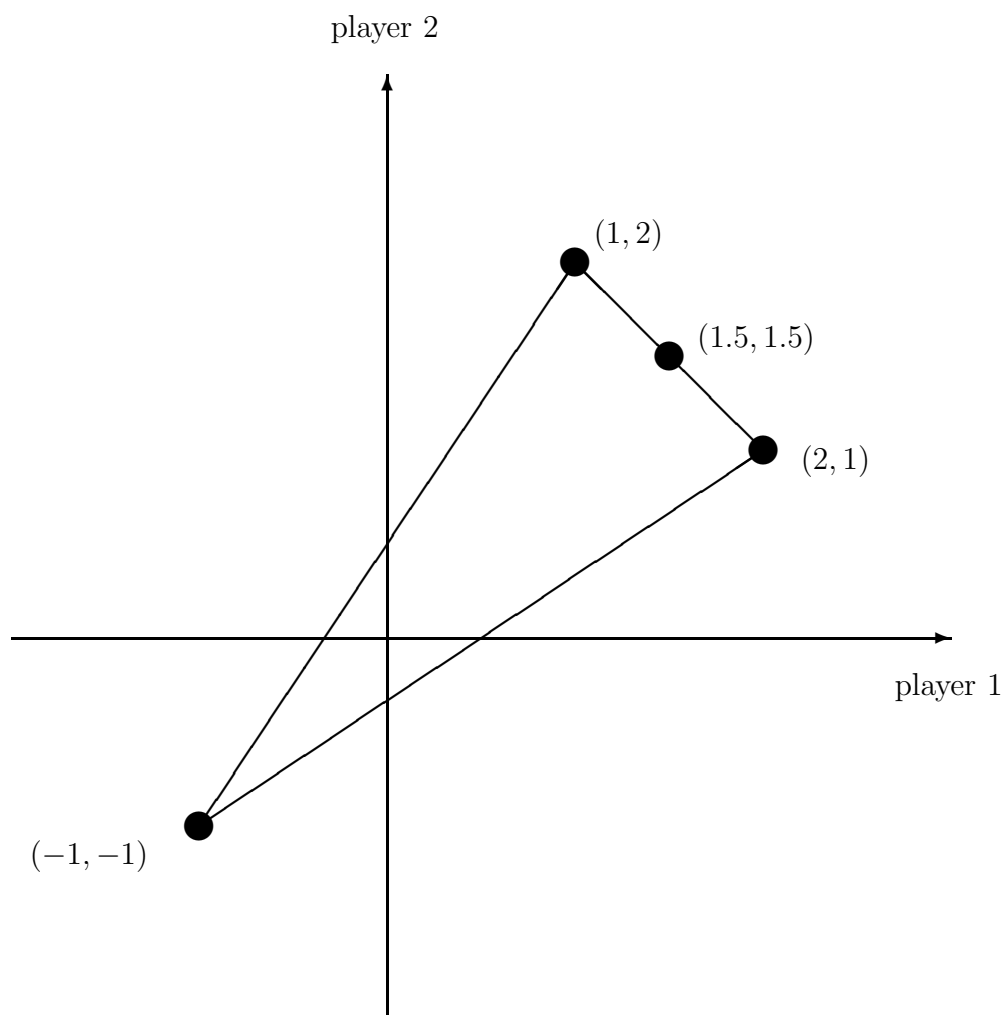
\begin{math}
\setlength{\unitlength}{2.5cm}
\thicklines
\begin{picture}(5,5.5)
\put(2,0){\vector(0,1){5}}
\put(1.7,5.2){$\mbox{player~2}$}
\put(0,2){\vector(1,0){5}}
\put(4.7,1.7){$\mbox{player~1}$}
\put(0.1,0.8){$( - 1, - 1 )$}
\put(1,1){\circle*{0.15}}
\put(1,1){\line(2,3){2}}
\put(3.1,4.1){$( 1, 2 )$}
\put(3,4){\circle*{0.15}}
\put(1,1){\line(3,2){3}}
\put(4.2,2.9){$( 2, 1 )$}
\put(4,3){\circle*{0.15}}
\put(4,3){\line(-1,1){1}}
\put(3.5,3.5){\circle*{0.15}}
\put(3.65,3.55){$( 1.5, 1.5 )$}
\end{picture}
\end{math} \\

\medskip
A similar extension of payoff outcomes to a closed convex set is possible through cooperation
in the case of arbitrary n-person mixed strategy games (2.5). Therefore, in this way they all
enter within the Nash bargaining or cooperation model. \\ \\

{\bf Remark 6} \\

In the particular case of a two person zero sum game in its mixed strategy form, the convex
set $H$ has the simple form

$$ H ~=~ \{~ ( x , -x ) ~|~ a ~\leq~ x ~\leq~ b ~\} $$

\medskip
for certain $a, b \in \mathbb{R},~ a ~\leq~ b$.

\hfill $\Box$ \\ \\

Let us now return to the Nash bargaining or cooperation model. According to the rules of the
game, let us assume that the players can secure the point $~h_0 ~=~ ( h_{\,0\,1},~ h_{\,0\,2} )
\in H~$ given in the next figure, In other words, completely independent of one another,
player 1 can make sure to receive at least $h_{\,0\,1}$, and similarly, player 2 can be sure
to receive at least $h_{\,0\,2}$. \\

The issue, therefore, is to what extent can the players improve on the outcome $~h_0~=~
( h_{\,0\,1},~ h_{\,0\,2} ) \in H~$ by cooperation ? \\

In this regard, Nash suggested {\it four axioms} presented below, which were supposed to
impose themselves upon all rational players due to their obviously fair and natural features,
and thus set up a corresponding cooperation between them. Let us therefore have a more
detailed look at what is happening, before we present the axioms. \\

\newpage

\bigskip
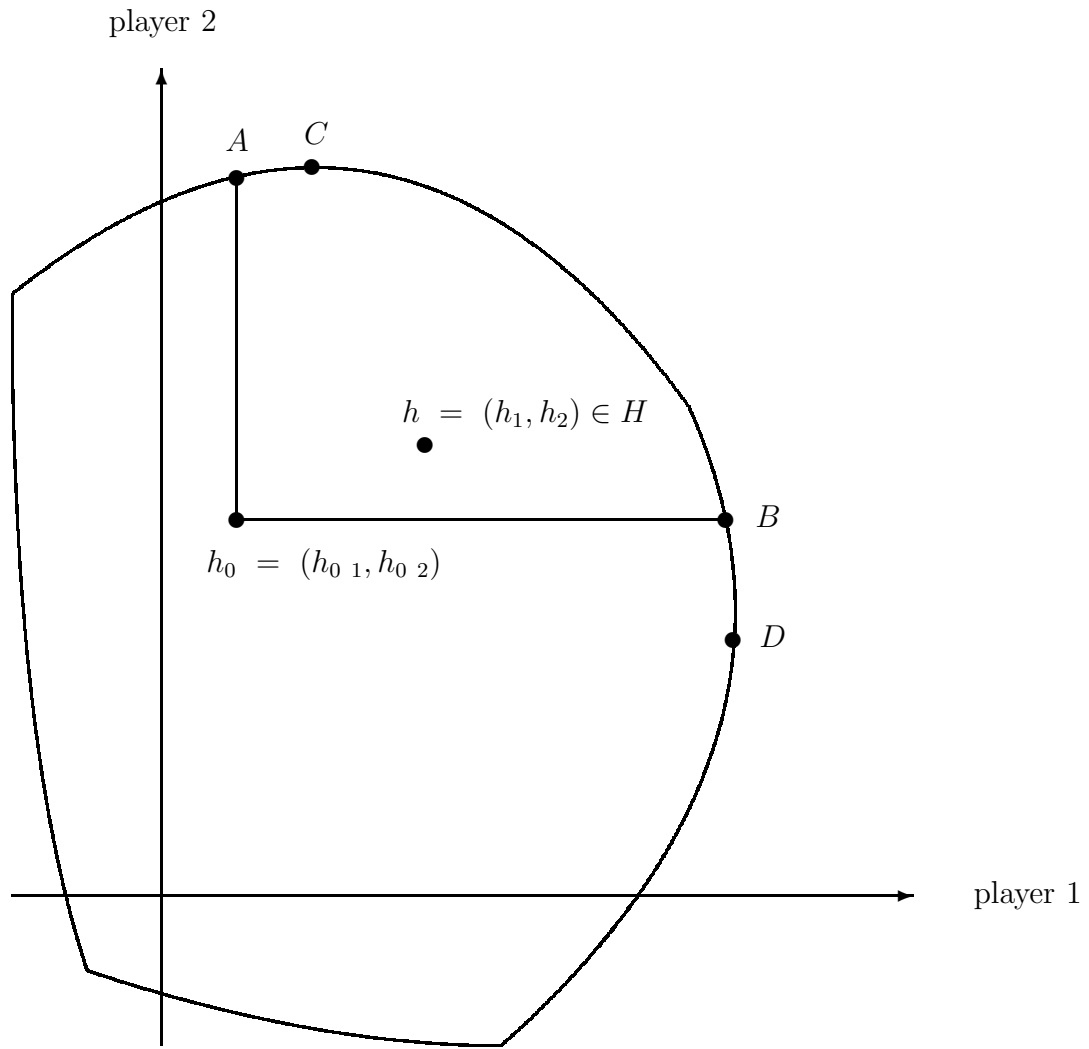
\begin{math}
\setlength{\unitlength}{1cm}
\thicklines
\begin{picture}(15,14)
\put(2,0){\vector(0,1){13}}
\put(1.3,13.5){$\mbox{player~2}$}
\put(0,2){\vector(1,0){12}}
\put(12.8,1.9){$\mbox{player~1}$}
\qbezier(0,10)(5,14)(9,8.5)
\qbezier(9,8.5)(11,4)(6.5,0)
\qbezier(6.5,0)(4,0)(1,1)
\qbezier(1,1)(0,4)(0,10)
\put(2.6,6.3){$h_0 ~=~ ( h_{0~1}, h_{0~2} )$}
\put(5.5,8){\circle*{0.2}}
\put(5.2,8.3){$h ~=~ ( h_1, h_2 ) \in H$}
\put(3,7){\circle*{0.2}}
\put(3,7){\line(1,0){6.5}}
\put(9.5,7){\circle*{0.2}}
\put(9.9,6.9){$B$}
\put(3,7){\line(0,1){4.6}}
\put(3,11.55){\circle*{0.2}}
\put(2.85,11.9){$A$}
\put(9.6,5.4){\circle*{0.2}}
\put(9.95,5.3){$D$}
\put(4,11.7){\circle*{0.2}}
\put(3.9,12){$C$}
\end{picture}
\end{math} \\ \\

\medskip
With $~h_0 ~=~ ( h_{\,0\,1},~ h_{\,0\,2} ) \in H~$ available even in the worst case, the
interest of the players is to cooperate in order to be able to move to a better point $~h ~=~
( h_1,~ h_2 ) \in H~$, that is, for which both $h_1 \geq h_{\,0\,1}$ and $h_2 \geq h_{\,0\,2}$
hold, and which we shall denote in short by \\

(3.1) \quad $ ~h_0 ~=~ ( h_{\,0\,1},~ h_{\,0\,2} ) ~\leq~ ~h ~=~ ( h_1,~ h_2 ) $ \\

Clearly, for the purpose of improving on $h_0$, only the subset $H^\prime$ of $H$ which is
contained within the lines $h_0, A, C, B, h_0$ is of interest. Furthermore, only the
intersection $H^*$ of this subset $H^\prime$ with the so called {\it Pareto maximal} subset of
$H$ is actually relevant. Here the Pareto maximal subset is given by the line $C, D$, thus
$H^*$ will result in the line $C, B$. Indeed, for any other $~h ~=~ ( h_1, h_2 ) \in
H^\prime~$ there exists $h^* ~=~ ( h^*_1, h^*_2 ) \in H^*$ such that $h^*_1 \geq h_1$ and
$h^*_2 \geq h_2$. \\

For convenience, let us recall that the Pareto maximal subset of any set H in $\mathbb{R}^2$
is the set of {\it nondominated} points $h = ( h_1, h_2 )$ in $H$, that is, those which
satisfy the condition \\

(3.2) \quad $ h ~\leq~ h^\prime \in H ~~~\Longrightarrow~~~ h^\prime ~=~ h $ \\

Here we should recall that, in general, reducing the initial convex, closed and bounded set
$H$ of payoffs to its Pareto maximal subset is the concept of solution in a two person nonzero
sum game which was suggested by von Neumann \& Morgenstern. \\
In case there is a guaranteed payoff $ ~h_0 ~=~ ( h_{\,0\,1},~ h_{\,0\,2} ) \in H$, then as
above, the Pareto maximal subset of $H$ can further be reduced to what is denoted by $H^*$. \\

Such a solution concept, however, is not satisfactory due to the following two facts :

\begin{enumerate}

\item In general, the subset $H^*$ is not a single point, and it can contain infinitely many
points.

\item When the two players are reduced to choosing one point in $H^*$, they face a situation
in which there {\it cannot} be any cooperation. Indeed, given any two different points $a =
( a_1, a_2),~ b = ( b_1, b_2 ) \in H^*,~ a \neq b$, we must have either \\

$~~~~~~ a_1 ~<~ b_1 ~~\mbox{and}~~ a_2 ~>~ b_2 $ \\

or \\

$~~~~~~ a_1 ~>~ b_1 ~~\mbox{and}~~ a_2 ~<~ b_2 $ \\

Thus always one of the two players must end up worse by going from $a$ to $b$, or from $b$ to
$a$.

\end{enumerate}

\medskip
{\bf Remark 7} \\

In the case of two person zero sum games, see Remark 6, we have $H^* ~=~ H$, that is, the
whole of original set $H$ of payoff outcomes is the Pareto maximal set. \\
Fortunately, in this particular case, the celebrated 1928 Min-Max theorem of von Neumann can
always yield a unique solution. However, such games do not allow any cooperation, therefore
they are outside of the framework of the Nash bargaining or cooperation.

\hfill $\Box$ \\ \\

The problem which arose in the late 1940s was to find, at least in the case of two person
nonzero sum cooperative games, a {\it solution concept} which would always deliver a {\it
unique} solution, and thus go beyond the solutions given by Pareto maximal sets. \\
And here the Nash axioms came to the fore in a rather obvious manner. According to them,
cooperation, or bargaining, is given by any function $F$ which depends on the {\it convex,
closed} and {\it bounded} sets $H$ and the respective guaranteed payoffs $~h_0 ~=~
( h_{\,0\,1},~ h_{\,0\,2} ) \in H~$, and always has a corresponding {\it unique} value given
by \\

$~~~~~~ h^* ~=~ ( h^*_1, h^*_2 ) ~=~ F ( H, h_0 ) \in H $ \\

a value called the {\it Nash bargaining or cooperation solution}, provided that the following
{\it four axioms} are satisfied : \\

{\bf Axiom of Pareto Maximality}

$$ h_0 ~~\leq~~ F ( H, h_0 ) ~=~ h^* ~=~ ( h^*_1, h^*_2 ) \in H^* $$

\medskip
where above inequality means, see (3.1), the two inequalities, $h_{\,0\,1} \leq h^*_1$ and
$h_{\,0\,2} \leq h^*_2$. \\

{\bf Axiom of Symmetry}

$$ F ( sym~ H, sym~ h_0 ) ~=~ sym~ F ( H, h_0 ) $$

\medskip
where $sym : \mathbb{R}^2 ~\longrightarrow~ \mathbb{R}^2$ is defined by $sym ( x, y ) ~=~ ( y, x )$. \\

{\bf Axiom of Independence of Irrelevant Alternatives} \\

Given two convex subsets $G ~\subseteq~ H ~\subseteq~ \mathbb{R}^2$, then

$$ h_0,~ F ( H, h_0 ) \in G ~~~\Longrightarrow~~~ F ( G, h_0 ) ~=~ F ( H, h_0 ) $$

In other words, if both the guaranteed payoff $h_0$, and the Nash cooperative or bargaining
solution $F ( H, h_0 )$ happen to belong to a smaller convex set $G ~\subseteq~ \mathbb{R}^2$,
then the solution $F ( G, h_0 )$ for that smaller set remains the same with the original
solution $F ( H, h_0 )$ for the larger set $H$. Or equivalently, if from the original set $H$
we leave out parts which do not contain the guaranteed solution $h_0$ and the solution
$F ( H, h_0 )$, then such an elimination does not affect the solution. \\

Finally, an axiom about the independence of changes of units in measuring the payoffs. \\

{\bf Axiom of Scaling Invariance} \\

Given any function $scale : \mathbb{R}^2 ~\longrightarrow~ \mathbb{R}^2$ such that $scale
( x, y ) ~=~ ( a + b x, c + d y )$, with $a, b, c, d \in \mathbb{R}$ and $b, d > 0$, we have

$$ F ( scale~ H, scale~ h_0 ) ~=~ scale~ F ( H, h_0 ) $$

\medskip
The remarkable fact about the above four Nash axioms is presented in \\ \\

{\bf Theorem ( Nash)} \\

There exists a {\it unique} function \\

$~~~~~~ ( H, h_0 ) ~\longmapsto~ F ( H, h_0 ) ~=~ h^* ~=~ ( h^*_1, h^*_2 ) \in H $ \\

which satisfies the four Nash axioms, and thus it always gives the unique Nash bargaining or
cooperation solution. \\

Equivalently, the Nash bargaining or cooperation solution $h^* ~=~ ( h^*_1, h^*_2 )$ is the
{\it unique} point in $H$ which satisfies the condition \\

(3.3) \quad $ ( h_1 - h_{0~1} ) ( h_2 - h_{0~2} ) ~\leq~
                           ( h^*_1 - h_{0~1} ) ( h^*_2 - h_{0~2} ) $ \\

for all $h ~=~ ( h_1, h_2 ) \in H,~ h ~\geq~ h_0$. \\

In other words we have the following {\it product maximization characterization} of the {\it
unique} Nash bargaining or cooperation solution $h^* ~=~ ( h^*_1, h^*_2 ) \in H$, namely \\

(3.4) \quad $ \mbox{max}~ ( h_1 - h_{0~1} ) ( h_2 - h_{0~2} ) ~=~
                             ( h^*_1 - h_{0~1} ) ( h^*_2 - h_{0~2} ) $ \\

where the maximum is taken over all $h ~=~ ( h_1, h_2 ) \in H,~ h ~\geq~ h_0$. \\ \\

{\bf Remark 8} \\

In Harsanyi there is an interpretation, going back to Zeuthen, of the above characterization
by product maximization of the unique Nash bargaining or cooperation solution. Namely, let us
assume that each of the two players $~i \in P ~=~ \{~ 1, 2 ~\}$ has chosen a respective point
$h_{*~i} ~=~ ( h_{*~i~1}, h_{*~i~2} )$ in the Pareto maximal set $H^*$. Further, let us assume
that $h_{*~1~1},~ h_{*~2~2} ~>~ 0$. Then in case \\

(3.5) \quad $ ( h_{*~1~1} - h_{*~2~1}) ~/~ h_{*~1~1} ~\leq~
                              ( h_{*~2~2} - h_{*~1~2} ) ~/~ h_{*~2~2} $ \\

it means that the {\it relative loss} to player 1, when giving up his or her choice $h_{*~1}$
and accepting instead the choice $h_{*~2}$ of player 2, is not larger than the {\it relative
loss} to player 2 when switching from $h_{*~2}$ to $h_{*~1}$. Thus in this case it is fair to
expect player 1 to make the switch either to the choice of player 2, or to a third choice
which may give a {\it smaller} relative loss to player 2. In case the opposite inequality
holds between the respective relative losses, then player 2 is expected to switch. It follows
that a choice which is acceptable to both players may emerge after an iteration which brings
the above inequality and it opposite as near as possible to equality. \\

However, we can note that the above inequality (3.5) is equivalent to the inequality \\

(3.6) \quad $ h_{*~1~1}~ h_{*~1~2} ~\leq~ h_{*~2~1}~ h_{*~2~2} $ \\

Now as noted, bargaining need not mean accepting the other player's offer. Instead, one can
suggest a third point $h_{*~3} ~=~ ( h_{*~3~1}, h_{*~3~2} )$ in $H^*$ such that the
corresponding product  $h_{*~3~1}~ h_{*~3~2}$ is not less than any of the the two products in
(3.6). Thus the corresponding bargaining iteration process can be seen as {\it maximizing}
such products, as happens in (3.4) in the Nash bargaining or cooperation model. \\ \\

{\bf Remark 9} \\

There have been a number of criticisms of the four Nash axioms, and in particular, of the
axiom of independence of irrelevant alternatives, and also to some extent of the axiom of
scaling invariance, as well as of symmetry, see Luce \& Raiffa [pp. 128-134]. \\
Here for the sake brevity, we recall in short only one of the criticisms of the axiom of
independence of irrelevant alternatives. \\
Let us consider the following pair of two person nonzero sum cooperative games which are given
respectively by the closed and bounded convex sets

\bigskip
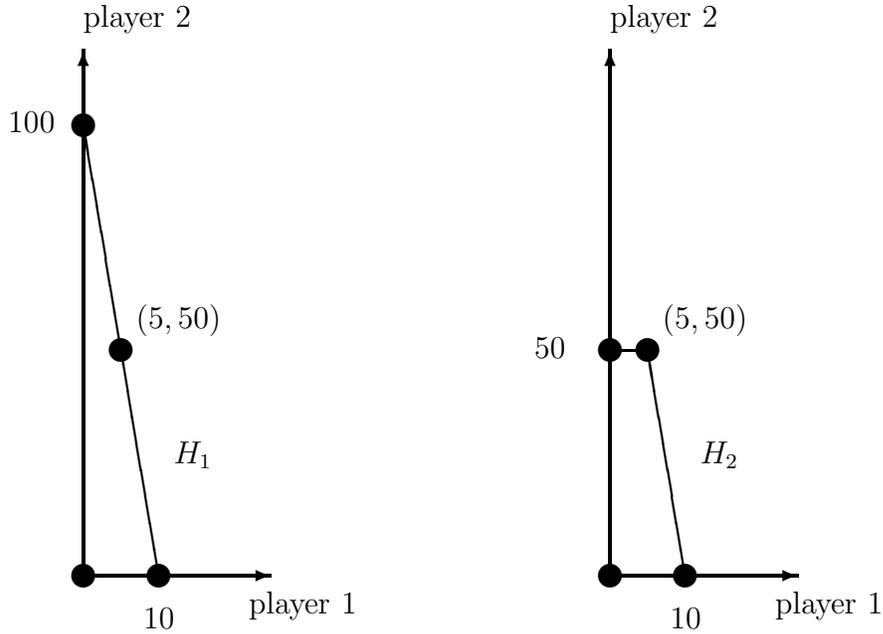
\begin{math}
\setlength{\unitlength}{1cm}
\thicklines
\begin{picture}(15,8)
\put(2,0){\vector(0,1){7}}
\put(2,7.3){$\mbox{player~2}$}
\put(2,0){\vector(1,0){2.5}}
\put(4.2,-0.5){$\mbox{player~1}$}
\put(3,0){\line(-1,6){1}}
\put(2,0){\circle*{0.3}}
\put(1,5.9){$100$}
\put(2,6){\circle*{0.3}}
\put(3,0){\circle*{0.3}}
\put(2.8,-0.7){$10$}
\put(3.2,1.5){$H_1$}
\put(2.5,3){\circle*{0.3}}
\put(2.7,3.3){$(5,50)$}
\put(9,0){\vector(0,1){7}}
\put(9,7.3){$\mbox{player~2}$}
\put(9,0){\vector(1,0){2.5}}
\put(11.2,-0.5){$\mbox{player~1}$}
\put(10,0){\line(-1,6){0.5}}
\put(9,0){\circle*{0.3}}
\put(8,2.9){$50$}
\put(9,3){\circle*{0.3}}
\put(9,3){\line(1,0){0.5}}
\put(9.5,3){\circle*{0.3}}
\put(9.7,3.3){$(5,50)$}
\put(10,0){\circle*{0.3}}
\put(9.8,-0.7){$10$}
\put(10.2,1.5){$H_2$}
\end{picture}
\end{math} \\ \\

\medskip
where the set $H_2$ is obtained from $H_1$ by cutting off the part above the horizontal line
at 50. In this case, both games have the {\it same} Nash bargaining or cooperation solution
given by \\

(3.7) \quad $ h^* ~=~ ( h^*_1, h^*_2 ) ~=~ ( 5, 50 ) $ \\

which obviously satisfies as well the axiom of independence of irrelevant alternatives. \\

However, this situation in (3.7) can bring about the following objections :

\begin{enumerate}

\item Player 2 can claim to receive more in game $H_1$ than in game $H_2$, since in game $H_1$
he or she could in principle receive twice as much as in game $H_2$.

\item Conversely, player 1 can complain that player 2 should receive less in game $H_2$ than
in game $H_1$, since in game $H_2$ the possibilities of player 2 were in principle diminished
no less than twice.

\end{enumerate}

And to add to the difficulties in accepting the four Nash axioms of bargaining or cooperation,
one should mention that even the axiom of Pareto maximality, which at first may appear to be
hardly at all controversial, does in fact lead to significant difficulties, see Luce \&
Raiffa [pp. 133, 134]. \\

In this way, it becomes once again obvious that even at conceptual levels, and even in the
simplest possible cases, such as those of cooperation which happen in two person games, there
are considerable difficulties involved. \\ \\

{\bf 4. Arrow's Impossibility Result} \\ \\

As seen above, it is not a trivial task to manage within a rational approach, such as game
theory for instance, the typically conflicting or competing objectives of a number of
autonomous players. And one of the difficulties may appear from the very beginning and at a
fundamental level, namely, with finding suitable {\it solution concepts} which take into
account under reasonable and general enough conditions the various individual objectives and
possibilities involved. \\
For instance, as it happens with the Nash equilibrium, the respective solution concept proves
to overturn the alleged non-cooperative nature of the games (2.1), (2.5). Also, in the Nash
bargaining or cooperation model, which only involves two players, what may appear to be
natural requirements about a solution do in fact lead to a variety of obvious difficulties. \\

In this section, we consider briefly the situation in what is usually called social,
collective or group choice. The awareness about the significant difficulty - even on the basic
conceptual levels - in setting up a reasonable aggregation of often conflicting or competing
objectives of autonomous agents has a longer history. \\
As mentioned, one of its cases, namely, the so called {\it Voter's Paradox}, dates back more
than two centuries, and has been addressed by persons like Condorcet, Laplace and Lewis
Carroll, among others, see Mirkin. One of the simplest forms of this paradox is as follows.
Three voters A, B and C have to rank in order of their own respective preferences three
alternatives a, b and c. The following situation can arise

$$ \begin{array}{l}
                  A ~~\mbox{has the preferences}~~ a < b < c \\
                  B ~~\mbox{has the preferences}~~ b < c < a \\
                  C ~~\mbox{has the preferences}~~ c < a < b
     \end{array} $$

\medskip
In such a case each of the following preferences $a < b$,~ $b < c$ and $c <a$ will be chosen
by 2 out of the 3 voters, thus by a majority of voters. And as often, in case the respective
preference relation $<$ is transitive and not reflexive, thus it is a {\it strict preference}
relation, then we obviously face the paradox that a majority considers $a < a$. By the same
token, a majority will also consider $b < b$ and $c < c$. \\

Needless to say, the appearance of a provable impossibility in any branch of science tends to
create a strong and long lasting impact. And the existence of a rigorous mathematical proof
for such an impossibility can only give further weight to such a result. Russell's paradox of
1901 in set theory, for instance, led to a long period of important studies in the foundations
of mathematics. The effects of G\"odel's 1931 impossibility result on formal systems in logic
has its reverberations still felt in mathematics. \\

Arrow's 1950 impossibility result has become one of the most important influences in group,
social or collective choice theory. Interestingly enough, Arrow's first formulation in 1950 of
that impossibility proved to contain certain errors, as pointed out by Blau. However, after
suitable corrections, a somewhat weaker form of impossibility could be reestablished, see
Arrow [2]. And in fact, subsequently, a whole range of related impossibility results has been
developed, see Kelly. We present here in a simple form the basic impossibility result of Arrow,
see Vincke [1], or Luce \& Raiffa. \\

Let $A$ be a finite set with at least three elements which represents all the possible choices.
We shall consider on $A$ binary relations $\leq$ which are transitive and complete, the latter
meaning that for every $a, b \in A$, we have at least one of the relations $a \leq b$, $ a =
b$, or $b \leq a$. For such a binary relation $\leq$, we shall define the associated strict
binary relation $<$, where for $a, b \in A$ we have $a < b$, if and only if $a \leq b$ and at
the same time we do not also have $b \leq a$. Further, we denote by $P_A$ the set of all such
binary relations on $A$. \\

Let us assume now that there are $n \geq$ 3 autonomous individuals $i$, each of whom can
freely choose a binary relation $\leq_i$ from $P_A$, which expresses his or her preferences in
the sense that, for $a, b \in A$, the individual $i$ prefers $b$ to $a$, if and only if $a
\leq b$. In this way, we are interested in the set \\

(4.1) \quad $ P_A^n ~=~ P_A \times ~.~.~.~ \times P_A $ \\

where the Cartesian product contains $n$ factors. Indeed, every n-tuple $( \leq_1, ~.~.~.~ ,
\leq_n ) \in P_A^n$ will be nothing else but the expression of one possible case of the
preferences of all the $n$ individuals involved. \\

Now, an {\it aggregation procedure} is by definition any function \\

(4.2) \quad $ f : P_A^n ~\longrightarrow~ P_A $ \\

which associates to every n-tuple $( \leq_1, ~.~.~.~ , \leq_n ) \in P_A^n$ expressing the
preferences of the $n$ individuals, a group preference

$$ \leq ~=~ f ( \leq_1, ~.~.~.~ , \leq_n ) \in P_A $$

\medskip
on $A$. Clearly, the problem is that, typically, the $n$ individual preferences

$$ \leq_1, ~.~.~.~ , \leq_n $$

\medskip
may conflict, and therefore, their aggregation into a group preference

$$ \leq ~=~ f ( \leq_1, ~.~.~.~ , \leq_n ) \in P_A $$

\medskip
is not a trivial task, provided that such an aggregation intends to take into account as much
as possible all the $n$ individual preferences. And the natural way to do so is to require
that the aggregation procedures $f$ to be used satisfy certain conditions which are seen to be
fair. \\

The conditions Arrow considered in this regard are the following three. First is the \\

{\bf Unanimity Condition} \\

(4.3) \quad $ \begin{array}{l}
                \forall~~~  ( \leq_1, ~.~.~.~ , \leq_n ) \in P_A^n,~~ a, b \in A ~: \\ \\
                ~~~~~ (~~ a <_i b ~~\mbox{for}~ 1 \leq i \leq n ~~~ )
                               ~~~\Longrightarrow~~~ a < b
              \end{array} $ \\

\medskip
where $\leq ~=~ f ( \leq_1, ~.~.~.~ , \leq_n )$. \\

Second is the \\

{\bf Condition of Independence of Irrelevant Alternatives} \\

(4.4) \quad $ \begin{array}{l}
                  \forall~~~ ( \leq_1, ~.~.~.~ , \leq_n ),~
                                ( \leq_1^\prime, ~.~.~.~ , \leq_n^\prime )
                                  \in P_A^n,~~ a, b \in A ~: \\ \\
                 ~~~~ (~~~ ( \leq_1, ~.~.~.~ , \leq_n ) ~=~
                 ( \leq_1^\prime, ~.~.~.~ , \leq_n^\prime ) ~~\mbox{on}~
                                          \{~ a,~ b ~\} ~~~) ~~~\Longrightarrow~~~  \\ \\
     ~~~~~~~~~~~~~~~~\Longrightarrow~~~  (~~~ \leq ~=~ \leq^\prime ~~\mbox{on}~
                                                               \{~ a,~ b ~\} ~~~)
              \end{array} $ \\

\medskip
where $\leq ~=~ f ( \leq_1, ~.~.~.~ , \leq_n )$ and $\leq^\prime ~=~
f ( \leq_1^\prime, ~.~.~.~ , \leq_n^\prime )$. \\

Finally, as the third one comes the \\

{\bf Condition on the Inexistence of a Dictator} \\

which means that is {\it none} of the individuals $i^{\#}$, with $1 \leq i^{\#} \leq n$ enjoys
the property \\

(4.5) \quad $ \begin{array}{l}
               \forall~~~  ( \leq_1, ~.~.~.~ , \leq_n ) \in P_A^n,~~ a, b \in A ~: \\ \\
                             ~~~~ a <_{i^{\#}} b ~~~\Longrightarrow~~ a < b
              \end{array} $ \\

\medskip
where $\leq ~=~ f ( \leq_1, ~.~.~.~ , \leq_n )$. \\

What follows is the celebrated {\it impossibility} result of Arrow, see Vincke [1], Luce \&
Raiffa, Kelly, or Arrow [2] \\ \\

{\bf Theorem ( Arrow ) } \\

Suppose that the set $A$ of possible choices has at least three elements, and $n \geq$ 3, that
is, there are at least three individuals. \\

Then there do {\it not} exist aggregation procedures $f : P_A^n ~\longrightarrow~ P_A$ in
(4.2) which satisfy the Unanimity and Independence conditions, and are without a
Dictator. \\ \\

{\bf Remark 10} \\

It is important, and also quite easy, to see that the crux of Arrow's impossibility is in the
fact that the preferences $\leq$ in $P_A$ considered on the set $A$ are assumed to be not only
complete, but also {\it transitive}. Indeed, let us consider on $A$ the larger set $Q_A$ of
binary relations on $A$ which need no longer be transitive, but only complete. Further,
instead of the original problem of finding an aggregation function (4.2), let us consider the
more general one of finding an aggregation function \\

(4.6) \quad $ g : Q_A^n ~\longrightarrow~ Q_A $ \\

Now, in view of the fact that the sought after aggregated preference $ \leq ~=~
g ( \leq_1, ~.~.~.~ , \leq_n ) \in Q_A$ need {\it no} longer be transitive, this aggregation
can be done very easily, for instance, according to {\it majority} rule, like in the example
of Condorcet. Namely, for $a, b \in A$, we can simply define \\

(4.7) \quad $ a \leq b ~~~\Leftrightarrow~~~
                  \mbox{car}~ \{~ i ~~|~~ 1 \leq i \leq n,~ a \leq_i b ~\} \geq~ n / 2 $ \\

Then it follows easily that, indeed, $\leq~ \in O_A$, and furthermore, $\leq$ satisfies all
the three conditions (4.3) - (4.5) in Arrow's theorem, namely, the Unanimity Condition, the
Condition of Independence of Irrelevant Alternatives, and for $n \geq$ 3, also the Condition
on the Inexistence of a Dictator. \\ \\

{\bf Remark 11} \\

It should be noted that in practical situations one can easily encounter individual
preferences which are {\it not} transitive. This can happen, for instance, when such
individual preferences result from multiple and conflicting objectives, see Rosinger [3,4] and
the literature cited there. \\

{\bf 5. Further on Cooperation} \\ \\

In section 3, some of the main issues related to cooperation, its possibilities and advantages
were considered within game theory. Here, in a similar manner, and in the case of social,
collective or group choice, several issues related to cooperation, its possibilities and
advantages are approached. \\
Clearly, one rather obvious rational way to overcome Arrow's impossibility, and do so without
a recourse to a dictator, is by a cooperation among the $n$ autonomous individuals involved,
cooperation aimed to produce certain modifications of the preferences of some of them. Of
course, such modifications may be compensated by suitable side payments, just as it happens in
games. \\

First, however, let us note that the number of possibilities for such modifications in
individual preferences - and thus for cooperation - are indeed {\it considerable}, since what
is involved in such modifications is to move from the original set \\

(5.1) \quad $ ( \leq_1, ~.~.~.~ , \leq_n ) \in P_A^n $ \\

of $n$ preferences of the respective individuals, to some other, this time more or less {\it
commonly} agreed upon new set of $n$ individual preferences \\

(5.2) \quad $ ( \leq_1^\prime, ~.~.~.~ , \leq_n^\prime ) \in P_A^n $ \\

Indeed, when the move is made from the original individual preferences in (5.1) to those in
(5.2), not all $n$ individuals need to effect modifications. Also, the modifications
effectuated do not necessarily need the agreement of everybody. Certainly, some such
modifications of preferences can be made voluntarily, either by certain single individuals or
by a number of them. And the motivation for such a modification in one's individual
preferences can range over a large spectrum, possibly including also the intent to avoid the
introduction of a dictator. \\

In this way, the {\it existence} of the possibility for such a move in preferences - and thus
for cooperation - is clearly there. Moreover, unlike in games, it is also very clear, as well
as simple {\it what} the $n$ individuals involved are supposed to do. Of course, as to {\it
how} they are to effectuate such a move in preferences, this remains just about as a complex
an issue, as it is in game theory. \\

Let us note here as a further argument supporting the possibility of cooperation, the {\it
extremely large number} of possibilities for a move in individual preferences. Indeed, let us
denote by $N$ the number of elements in $A$, that is, the number of possible choices which
each of the $n$ autonomous individuals have. Then it is obvious that \\

(5.3) \quad $ \mbox{car}~ P_A ~\geq~ N~ ! $ \\

hence \\

(5.4) \quad $ \mbox{car}~ P_A^n ~\geq~ ( N~ !~ )^n $ \\

Thus even in the case of, for instance, only 5 choices and 5 individuals, the size of $P_A^n$
is already {\it larger} than $10^{10}$. \\

It follows that it is not particularly rational to neglect, disregard or reject prior to the
aggregation of preferences, the possibilities for cooperation given by suitable modifications
of the original individual preferences in (5.1). \\ \\

{\bf 6. A Single Decision Maker with Multiple Conflicting Objectives} \\ \\

In view of the difficulties in games or social, collective and group choices to manage
conflict or competition, it may at first appear that, as mentioned, the process of making
decisions by one single decision maker may be so much more easy. And at first sight it may
also appear that Arrow's impossibility result would point in this direction. \\
In fact, as also mentioned, the case of one single decision maker, whom we shall denote by SDM,
can be seen as a game with one single player who plays against Nature. Consequently, he or she
can in a certain sense function as a dictator to himself or herself. \\

What actually happens in such a situation is often contrary to such a first perception. Namely,
the whole range of complexities and difficulties related to conflict, competition, cooperation,
bargaining, side payments, and so on, simply comes down as if crashing on the head of that SDM,
and thus become his or her internalized problems. \\

The one {\it major difference} is that in the case SDM behaves rationally, the {\it intention
to cooperate} is clearly there always, totally and unconditionally, since he or she is
cooperating this time with no one else but himself or herself. In this way, the case of such a
SDM brings with it the most {\it auspicious} context for cooperation. And then, the only issue
left, one that is not at all simple, is {\it how} to accomplish cooperation, which now becomes
but an issue of {\it competence}. \\

And as it turns out, this competence requires an understanding of the nature of the conflicts
involved in typical situations with multiple objectives, an understanding which can lead along
less than usual, and also somewhat counterintuitive directions of thinking and acting. \\

In this regard, two facts come to the fore from the beginning in typical situations with
multiple conflicting objectives, see Rosinger [3,4] :

\begin{itemize}

\item {\bf Fact 1.} There is {\it no}, and there {\it cannot} be a unique natural canonical
candidate for the concept of solution. And in fact, the very issue of choosing a solution
concept leads to a {\it meta-decision} problem which itself has multiple conflicting
objectives. \\

\item {\bf Fact 2.} The information contained in the preference structures involved, relative
to all other possible, and for instance, non-preference type information present in the
situation, tends {\it exponentially} to zero, as the number of conflicting objectives
increases. This phenomenon, which in fact is of a very simple nature related to the geometry
of {\it higher dimensional} Euclidean spaces, can be called the {\it Principle of Increasing
Irrelevance of Preference Type Information}, or in short PIIPTI, and it will be proved in
Appendix 1 at the end of this section. Its explicit identification was first presented in
Rosinger [4]. \\
A further complication comes from the frequent phenomenon that preferences, when they may be
found and expressed, tend {\it not} to be transitive.

\end{itemize}

What comes after the understanding of these two facts can be seen as being, so far, a kind of
{\it ultimate exercise in cooperation}, that is, in solving the conflicts involved, and doing
so within the context of the rational behaviour assumed to hold on the part of the SDM. \\

Let us address briefly these two facts which, so far, do not seem to be widely enough known. \\

When trying to solve the given decision problem with multiple conflicting objectives, the SDM
may at first have to set up and solve a meta-decision problem, namely, to choose an
appropriate {\it solution concept} for the initially given decision problem. \\
However, so often, the need for this first meta-decision stage is not realized. And then,
instead, one or another pet-solution concept which happens to be familiar is applied to the
decision problem, without any due consideration whether that particular solution concept may,
or may not be appropriate. \\

In order to illustrate in some more detail what is involved, let us consider the following
large and practically important class of decision making situations, when the SDM has to deal
with $n ~\geq~ 2$ typically conflicting objectives given by the utility functions \\

(6.1) \quad $ f_1, ~.~.~.~ , f_n ~:~ A ~~\longrightarrow~~ \mathbb{R} $ \\

and his or her aim is to maximize all of them, taking into account that most often such a
thing is not possible simultaneously, due to the conflicts involved. \\

Here, as before, the set $A$ describes the available choices, this time those which the SDM
has, and in this case $A$ may as well be an infinite set, for instance, some open or closed
bounded domain in a finite dimensional Euclidean space. \\
Clearly, the functions $f_i$ in (6.1) can be seen as utility functions, von Neumann \&
Morgenstern, Luce \& Raiffa, and as such, they generate preference relations on the set of
choices $A$. Namely, the preference relation $\leq_i$ corresponding to the utility function
$f_i$ is defined for $a, b \in A$, by \\

(6.2) \quad $ a \leq_i b ~~~\Leftrightarrow~~~ f_i ( a ) \leq f_i ( b ) $ \\

Needless to say, the situation described by (6.1) is not the most general one, since it is
possible to encounter cases when the objectives are not given by utility functions, or simply,
are not even quantifiable. However, the model in (6.1) can nevertheless offer an edifying
enough situation, in order to be able to obtain relevant insights into the nature and extent
of the complexities and difficulties which a SDM can face. Furthermore, it can also lead to
general enough solution methods, including ways to choose solution concepts, see Rosinger
[3,4]. \\

But to better illustrate the conundrum imposed upon a SDM, the following age old analogy may
be helpful in trying to describe his or her problem :

\begin{quote}

"How to run after, and hunt down each of two or more rabbits, and do so with one single gun
with one single bullet ?"

\end{quote}

So much for the ... joys ... in the life of Arrow's would be dictator tuned here into a
SDM ... \\

And now in view of the two facts mentioned above, let us indicate some of the more important
features involved which face a SDM. \\

First we note that, as we mentioned, the very {\it concept of solution} is not at all clear.
Indeed, there have been {\it three traditional} ways of dealing with the optimization of
multiple objectives. Also, since the 1970s, a number of other methods have been suggested.
Often however, such methods are not sufficiently compared with the already existing ones, and
instead, appear to be suggested as having a rather universal validity. Furthermore, such
methods do not sufficiently take into account the above mentioned two facts. \\

Let us now recall in short the three traditional methods in dealing with conflicting
objectives. \\

One of them is based on {\it priorities}. This method lists the objectives in a certain order,
and then tries to fulfill them as follows. First, the first objective is maximized, then
within the resources left, the second objective is maximized, and so on. Clearly, such an
approach is based on an up-front, instant, total, and often brutal, biased, or insufficiently
well considered elimination of conflict. And as such, it can easily lead to most of the
objectives, except for the first, and perhaps the very few next ones, being left without any
appropriate consideration. \\

The second method, aimed to improve on that of priorities, introduces an additional, overall
objective, which is constructed by a certain {\it weighted sum} of the initial conflicting
objectives, thus it is supposed to take into account all of them. In this way, the conflict is
again quite up-front and totally eliminated, and the problem is reduced to the optimization of
one single objective, a problem which is rather trivial in comparison. \\

The disadvantages of these two traditional methods of conflict resolution have been repeatedly
experienced in a large variety of practical situations. Chief among these disadvantages is the
fact that both methods lead to a far {\it too early, simple and total elimination} of
conflict. \\

As a third attempt, aimed to go beyond such disadvantages was suggested by Vilfredo Pareto,
who introduced the concept of {\it non-dominated solution}, see (3.2), which nowadays is also
called a {\it Pareto maximal} solution. This solution concept, with the help of suitable
additional analysis, such as for instance, trade-offs between various marginal utilities, was
supposed to lead to optimal decisions. \\

A deficiency of the Pareto solution concept is that, contrary to the method of priorities or
weighted objectives, it still leaves there, and not properly dealt with, {\it too much} of the
initial conflict, especially when $n \geq$ 3, where $n$ is the number of conflicting
objectives involved, see (6.1). \\
This is related to a property of finite dimensional Euclidean spaces with larger dimensions,
and it comes from the elementary, even if less widely known geometric fact that spheres in
such spaces have most of their {\it volume} concentrated in a thin layer near to their surface,
while the whole rest of the interior of such spheres only contains a small part of their
volume. In fact, this increase in the relative volume of the thin surface layer, when compared
with the whole of the rest of the volume, is {\it exponential} in the dimension of the
respective sphere, see Appendix 2 at the end of this section. \\

In other words

\begin{quote}

It is not worth buying ... higher dimensional water melons ... \\
And certainly, one should not peel ... higher dimensional fruit ...

\end{quote}

There is also a second, and yet simpler geometric argument leading to the deficiency of the
Pareto solution concept when the number $n$ of conflicting objectives is larger. This argument,
as seen below, concerns the relatively {\it high dimension} of the Pareto maximal, or
non-dominated sets. \\

Let us show now that within the model (6.1) of multiple conflicting objectives, when $n \geq$
3, the use of the Pareto maximal, or non-dominated sets as a solution concept does indeed
leave too much of the conflict unsolved, as far as the reduction in {\it dimensions} is
concerned. Indeed, let us denote by \\

(6.3) \quad $ B ~=~ \{~ ( f_1 ( a ), ~.~.~.~ , f_n ( a ) ) ~|~ a \in A ~\} ~\subseteq~
                                                                           \mathbb{R}^n $ \\

the set of n-tuples of outcomes $( f_1 ( a ), ~.~.~.~ , f_n ( a ) ) \in \mathbb{R}^n$ which
correspond to various choices $a \in A$ which the SDM can make. Typically, this is a bounded
and closed subset in $\mathbb{R}^n$, similar to what happens, when $n = 2$, in the case of two
person nonzero sum non-cooperative games, see section 3. \\

Let us now consider the Pareto maximal, or non-dominated subset of $B$, which will be given
by \\

(6.4) \quad $ B^P ~=~ \{~ b \in B ~~|~~ \forall~~ c \in B ~:~~ b ~\leq~ c
                                            ~~\Longrightarrow~~ c ~=~ b ~\} $ \\

where, as in (3.1), here correspondingly in general, for $b = ( b_1, ~.~.~.~ , b_n ),~ c =
( c_1, ~.~.~.~ , c_n ) \in \mathbb{R}^n$, we denote \\

(6.5) \quad $ b \leq c $ \\

if and only if $b_i \leq c_i$, with $1 \leq i \leq n$. \\

Now, it is often the case that the set $B$ of outcomes has a nonvoid interior, and thus it is
an n-dimensional subset of $\mathbb{R}^n$. Consequently, the Pareto maximal subset $B^P$ is an
(n-1)-dimensional subset of $B$. In this way, by using the solution concept of Pareto maximal
or non-dominated subsets, all what is done is to go from a set $B$ with $n$ dimensions to a
subset of it which still has $n - 1$ dimensions. And while in the case of $n = 2$, that is, of
only two conflicting objectives, this reduces the dimension by half, on the other hand, when
$n \geq 3$, the relative reduction in dimension is merely $1 ~/~ n$, thus it becomes less and
less relevant, as $n$ may increase. \\

In conclusion, the solution concept given by the Pareto maximal or non-dominated subset may be
useful in the case of $n = 2$ conflicting objectives, similar to what happened with the Nash
bargaining or cooperative solution, see section 3. However, as $n$ increases, starting with
$n = 3$, this solution concept leaves more and more of the conflict unresolved, since both the
relative {\it volume} and relative {\it dimension} of the Pareto maximal, or non-dominated
sets become higher, as $n$ may increase. \\

Let us return now again to the above mentioned Fact 1, according to which {\it prior} to
solving the multiple conflicting decision making problem in (6.1), the SMD may have to solve a
{\it meta-decision} problem which itself may have multiple conflicting objectives. Namely, the
SDM may have to decide upon a {\it solution concept} which he or she will use in solving the
initial problem in (6.1). In this regard, as mentioned in Rosinger [3,4] and the literature
cited there, several {\it meta-objectives} concerning possible solution concepts have been put
forward. Among them are :

\begin{enumerate}

\item Ease in finding the solution to (6.1).

\item Fidelity in modelling the conflict in (6.1) when $n \geq 2$.

\item Confidence in the relevance of the solution obtained.

\item Use not only of preference type information, but of other possible information, such as
for instance, indifference.

\end{enumerate}

Further, in order to satisfy such meta-objectives about solution concepts, since the late
1970s a number of {\it interactive} methods have been developed for the solution of (6.1), see
Rosinger [3,4] and the literature cited there. \\

It is important to note here that, as argued next, such interactive methods lead to a new
class of so called {\it interactive solution concepts}, which are beyond the traditional ones
that were only giving what can be considered as being {\it a priori solution concepts}. \\
Indeed, in such interactive methods the corresponding solution concepts presented to the SDM
do not merely lead to input-output, black-box type algorithms, which the SDM can only use by
introducing the data about (6.1), and then having to wait passively until the result is
obtained at the output. \\
Instead, the SDM can during the execution of the algorithm interact with it, in order to
provide further information on, or make certain partial decisions about the conflicts involved
in (6.1). \\

In this way, one may add to the above meta-objectives the following one as well : \\

{~~~} 5. Ease, efficiency and effectiveness in the interaction process. \\

Details about such interactive solution concepts can be found in Rosinger [3,4], and the
literature cited there. \\

What is remarkable in the above, related to decision making with multiple conflicting
objectives, is that everything rational done in this respect can be seen as an instance of
{\it cooperation} aimed to deal properly with the conflicts involved. This includes the
meta-decision problem itself, which tries to identify a solution concept, as well as the
interactive type solution concepts. \\

On the other hand, since only one rational agent, that is, the SDM is involved, it may happen
that - due to the lack of presence of, and challenge from other rational agents - the SDM may
not always pursue consistently enough the ways of cooperation, in spite of the fact that he or
she is unreservedly cooperation minded.  \\
And the extent to which cooperation is possible has most likely not yet been figured out
satisfactorily. One indication in this regard can be given, for instance, by the lack of a
wide enough awareness about Fact 1 and Fact 2 mentioned at the beginning of this section. \\ \\

{\bf Appendix 1} \\

We start with a very simple geometric fact about {\it finite dimensional} Euclidean spaces
which can give a good insight into the more general result in (6.8). \\
On the n-dimensional Euclidean space $\mathbb{R}^n$, with $n \geq 1$, we consider the natural
partial order relation $\leq$ defined for elements $x = ( x_1, ~.~.~.~ , x_n ),~ y =
( y_1, ~.~.~.~ ,  y_n ) \in \mathbb{R}^n$, according to, see (6.5)

$$ x ~\leq~ y ~~~\Leftrightarrow~~ x_i ~\leq~ y_i,~~\mbox{with}~ 1 \leq i \leq n $$

\medskip
Let us denote by

$$ P_n ~=~ \{~ x  \in \mathbb{R}^n  ~~|~~ x ~\geq~ 0 ~\} $$

\medskip
the set of nonnegative elements in $\mathbb{R}^n$, corresponding to the partial order
$\leq$. \\

Then we can note that, for $n = 1$, the set $P_1$ is half of the space $\mathbb{R}^1 =
\mathbb{R}$. \\
Further, for $n = 2$, the set $P_2$ is a quarter of the space $\mathbb{R}^2$. \\
And in general, for $n \geq 1$, the set $P_n$ is $1 ~/~ 2^n$ of the space $\mathbb{R}^n$. \\

It follows that in an n-dimensional Euclidean space $\mathbb{R}^n$, if one is given an
arbitrary element $x  \in \mathbb{R}^n$, then the probability for this element $x$ to be
nonnegative is $1 ~/~ 2^n$, thus it tends {\it exponentially} to zero with $n$. \\
Consequently, the same happens with the probability that two arbitrary elements $x ,~ y \in
\mathbb{R}^n$ are in the relationship $x \leq y$, since this is obviously equivalent with the
condition $y - x \geq 0$. \\

This means that, when the number $n \geq 2$ of conflicting objectives increases, one can
expect a similar trend to happen with preference relations on the set of outcomes in (6.3),
which correspond to the $n$ conflicting objectives in (6.1), that phenomenon being at the root of
PIIPTI formulated in Fact 2 above. A respective result in this regard is presented in (6.8) in
the sequel. \\

Let us now assume that in (6.1) we have a finite set of choices, namely \\

(6.6) \quad $ A ~=~ \{~ a_1, ~.~.~.~ , a_m ~\},~~~ m ~\geq~ 2 $ \\

A natural single preference relation on $A$ corresponding to (6.1), and which may try to
synthesize the respective $n$ conflicting objectives, is given by a subset \\

(6.7) \quad $ S ~\subseteq~ A \times A $ \\

Here, for any $a,~ a^{\,\prime} \in A$, the SDM will prefer $a^{\,\prime}$ to $a$, in which
case we write $a ~\leq~ a^{\,\prime}$, or equivalently, $( a, a^{\,\prime} ) \in S$, if and
only if one has for each objective function $f_i$, with $1 \leq i \leq n$, either that $f_i
( a^{\,\prime} ) - f_i ( a ) ~>~ 0$ and it is not negligible, or ~$| f_i ( a^{\,\prime} ) -
f_i ( a ) |$~ is negligible. \\

One can note that the method of weighted sums also leads to such a single preference relation
$S$ for the SDM, while the method of priorities does in fact just about the same. Of course,
in these two latter cases the resulting preferences may be quite different from that in (6.7),
as well as from one another. \\

Obviously, there is a strong tendency to think that the SDM could solve the problem based on
such, or some similar, one single preference relation $S$ which he or she may be able to find,
such a tendency being not in the least a consequence of Arrow's impossibility result. And
needless to say that such a thinking comes mainly from the fact that there has not been a long
and sophisticated enough tradition, and therefore experience, in dealing with multiple
conflicting objectives. Instead, in such conflict situations there has traditionally been, as
mentioned, a tendency simply to eliminate the conflict too soon, and thus necessarily in
highly questionable ways. \\
And certainly, when $n = 1$, that is, when in (6.1) there is only one single objective, the
whole information is already given in one single corresponding preference relation $S$, and
thus it is enough to solve the problem. However, such a problem with $ n = 1$ is so simple
that it is {\it no} longer a problem of multiple conflicting objective decision making. \\

Let us therefore see more precisely how much information one single preference relation $S$
can carry, when the number $n$ of conflicting objectives in (6.1) becomes larger, and even if
only moderately so. This can be done quite easily by noting that in typical situations, we can
have the relation \\

(6.8) \quad $ \mbox{car}~ S ~/~ \mbox{car}~ ( A \times A ) ~=~ O ( 1 ~/~ 2^n ) $ \\

where for a finite set $E$ we denoted by "$\mbox{car}~ E$" the number of its elements. \\

The proof of (6.8) goes as follows, by using a combinatorial-probabilistic type argument. Let
us take any injective function $g : A \longrightarrow \mathbb{R}$, and denote by \\

(6.9) \quad $ S_g ~=~ \{~ ( a, a^{\,\prime} ) \in A \times A ~~|~~ g ( a ) ~\leq~
                                                              g ( a^{\,\prime} ) ~\} $ \\

which is its corresponding preference relation on $A$. Then obviously \\

(6.10) \quad $ \mbox{car}~ S_g ~=~ m ( m + 1 ) ~/~ 2 $ \\

Now given any subset $S \subseteq A \times A$, let us denote by $P ( S )$ the probability that
for an arbitrary pair $( a, a^{\,\prime} ) \in A \times A$, we have $( a, a^{\,\prime} ) \in
S$. Then clearly \\

(6.11) \quad $ P ( S_g ) ~=~ ( 1 + 1 ~/~ m ) ~/~ 2 $ \\

Let us assume now about the objective functions in (6.1) the following \\

(6.12) \quad $ f_1, ~.~.~.~ , f_n ~~\mbox{are injective} $ \\

and their corresponding sets of preferences \\

(6.13) \quad $ S_{f_1}, ~.~.~.~ , S_{f_n} ~~\mbox{are probabilisitically independent} $ \\

Then we obtain, see (6.11) \\

(6.14) \quad $ P  ( S_{f_1} \bigcap ~.~.~.~ \bigcap S_{f_n} ) ~=~
              P ( S_{f_1} ) ~.~.~.~ P ( S_{f_n} ) ~=~ ( 1 + 1 ~/~ m )^n ~/~ 2^n $ \\

And now (6.8) follows, provided that $n$ in (6.1) and $m$ in (6.6) are such that \\

(6.15) \quad $ ( 1 + 1 ~/~ m )^n ~=~ O ( 1 ) $ \\

which happens in many practical situations. \\

As for the independence condition (6.13), let us note the following. Let us assume that the
objectives $f_1$ and $f_2$ are such that for $a,~ a^{\,\prime} \in A$ we have \\

(6.16) \quad $ f_1 ( a ) ~<~ f_1 ( a^{\,\prime} ) ~~~\Leftrightarrow~~~
                                                f_2 ( a ) ~<~ f_2 ( a^{\,\prime} ) $ \\

then obviously $S_{f_1} = S_{f_2}$, hence (6.14) may fail. But clearly, (6.16) means that
$S_{f_1}$ and $S_{f_2}$ are {\it not} independent. In the opposite case, when \\

(6.17) \quad $ f_1 ( a ) ~<~ f_1 ( a^{\,\prime} ) ~~~\Leftrightarrow~~~
                                              f_2 ( a^{\,\prime} ) ~<~ f_2 ( a ) $ \\

then obviously \\

(6.18) \quad $ S_{f_1} \bigcap S_{f_2} ~=~ \{~ ( a, a ) ~|~ a \in A ~\} $ \\

and (6.14) may again fail. However (6.18) once more means that $S_{f_1}$ and $S_{f_2}$ are
{\it not} independent, since they are in total conflict with one another. \\ \\

{\bf Appendix 2} \\

For the sake of simplicity, let us assume that the set $B$ of outcomes in (6.3) is of the
form \\

(6.19) \quad $ B ~=~ \bigg \{~ b = ( b_1, ~.~.~.~ , b_n ) \in \mathbb{R}^n ~~
                     \bigg |~~ \begin{array}{l}
                          b_1, ~.~.~.~ , b_n ~\geq~ 0 \\
                          b_1 + ~.~.~.~ + b_n ~\leq~ L
                                \end{array} ~\bigg \} $ \\

for a certain $L > 0$. Then clearly the Pareto maximal subset of $B$ is \\

(6.20) \quad $ B^P ~=~ \{~ b = ( b_1, ~.~.~.~ , b_n ) \in \mathbb{R}^n ~~|~~
                                                     b_1 + ~.~.~.~ + b_n = L ~\} $ \\

Now for $0 < \epsilon < L$, the $\epsilon$-thin shell in $B$ corresponding to $B^P$ is given
by \\

(6.21) \quad $ B^P ( \epsilon ) ~=~ \bigg \{ b = ( b_1, ~.~.~.~ , b_n ) \in \mathbb{R}^n
           ~\bigg | \begin{array}{l}
   L - \epsilon ~\leq~  b_1 + ~.~.~.~ + b_n ~\leq~ \\
  ~~~~~~~~~~~~~~~~~~~~~~ \leq~ L
  \end{array} \bigg \} $ \\

A standard multivariate Calculus argument gives for the volume of $B$ in (6.19) the
relation \\

(6.22) \quad $ \mbox{vol}~ B ~=~ K_n~ L^n $ \\

where the constant $K_n > 0$, involving the Gamma function, does only depend on $n$ and not on
$L$ as well. In this way it is easy to see that \\

(6.23) \quad $ \mbox{vol}~ B^P ( \epsilon ) ~/~ \mbox{vol}~ B ~=~
                                        1 - ( 1 - \epsilon ~/~ L )^n $ \\

This shows, for instance, that in the 20-dimensional case, a shell with a thickness of only
5\% of the radius $L$ will nevertheless contain at least 63\% of the total volume. \\

In more simple and direct {\it geometric} terms the relation (6.23) means that :

\begin{quote}

{\bf ( VOL )}~~ "The volume of a multidimensional solid is  \\
                 \hspace*{2.2cm} mostly concentrated next to its surface."

\end{quote}

It may be instructive to note that relation (6.23) also has a {\it physical} interpretation,
as it explains the phenomenon of {\it temperature}, see Manin. Indeed, let us assume that a
certain simple gas has $n$ atoms of unit mass. Then their {\it kinetic energy} is given by \\

(6.24) \quad $ E ~=~ \Sigma_{1 \leq i \leq n}~ v_i^2 ~/~ 2 $ \\

where $v_i$, with $1 \leq i \leq n$, are the velocities of the respective atoms. Therefore,
for a given value of the kinetic energy $E$, the state of the gas is described by the vector
of $n$ velocities, namely \\

(6.25) \quad $ v ~=~ ( v_1, ~.~.~.~ , v_n ) \in S_n ( \sqrt ( 2 E ) ) $ \\

where for $L > 0$, we denoted by \\

(6.26) \quad $ S_n ( L ) ~=~ \{~ x = ( x_1, ~.~.~.~ , x_n ) \in \mathbb{R}^n ~~|~~
                                                   x_1^2 + ~.~.~.~ + x_n^2 = L^2 ~\} $ \\

the ( n - 1 )-dimensional surface of the n-dimensional ball with radius $L$ in
$\mathbb{R}^n$. \\

Now we can recall that in view of the Avogadro number, under normal conditions for a usual
macroscopic volume of gas, one can have \\

(6.27) \quad $ n ~>~ 10^{20} $ \\

Therefore, the above property {\bf ( VOL )} which follows from (6.23) is very much manifest.
Let us then assume that a small thermometer with a thermal energy $e$ negligible compared to
$E$ is placed in the gas. Then the state (6.25) of the gas will change to a new state \\

(6.28) \quad $ v ~=~ ( v_1, ~.~.~.~ , v_n ) \in S_n ( \sqrt ( 2 E^\prime ) ) $ \\

However, in view of property {\bf ( VOL )}, it will follow with a high probability that \\

(6.29) \quad $ E^\prime ~\approx~ E $ \\

And it is precisely this {\it stability} or {\it rigidity} property (6.29) which leads to the
phenomenon of {\it temperature} as a macroscopically observable quantity. \\ \\

{\bf 7. Conclusions} \\ \\

Situations involving the actions of conscious rational agents who aim to optimize certain
outcomes are approached in three mathematical theories, namely,
the theory of games, the theory of social, collective or group choice, and decision theory. In
games, there are two or more such conscious and rational agents, called players, who are
interacting according to the given rules. And except for that, they are free and independent,
and there is no overall authority who could influence in any way the players. In social,
collective or group choice, again, there are two or more conscious and rational agents with
their given individual preferences. Here however, the issue is to find a mutually acceptable
aggregation of those preferences. And such an aggregation is seen as being done by an outsider.
Finally, decision theory can be seen as a game in which a conscious rational agent plays
against Nature. \\

Cooperation is usually seen as pertaining rather exclusively to games. Due to the extreme
complexities involved, however, a systematic, deep and far reaching enough study of
cooperation has not yet been accomplished within game theory. \\

Although in ways not entirely identical with those in games, cooperation can take place as
well in situations of social, collective or group choice. And here, ever since the celebrated
1950 Arrow impossibility result, there is a strong motivation for implementing cooperation,
since the alternative in general is the introduction of a dictator. \\

One of the main sources of the extreme complexities related to cooperation in games and in
social, collective or group choice is in the presence of two or more conscious and rational
agents. And as seen in Binmore [1-3], such a presence in games can already take things well
beyond the reach of any algorithmic approach, even if cooperation is not among the main issues
pursued. Also, in social, collective or group choice, the only way to overcome the extreme
complexities created by the presence of three or more individual preferences is by the
introduction of a dictator, as proved in the mentioned result of Arrow. \\
Yet the presence of such extreme complexities related to cooperation need not lead to a
situation where one has in fact given up completely and for good on the study and
implementation of cooperation, especially if this giving up has happened rather by default. \\

In the case of one single decision maker, who in decision theory is seen as a conscious
rational player, playing alone against Nature, one may seem at first two have a situation
which enjoys all the advantages that are missing both in games, and in social, collective or
group choice. Indeed, it may at first appear that such a single decision maker does not have
to put up with one or more other autonomous agents. And also, as the single player, he or she
can automatically be seen as a dictator as well, since there is no other conscious agent out
there to protest and actively oppose, least of all what is called Nature in such a context. \\
It would, therefore, appear that in decision theory one has it very easy, and in particular,
one has no need to consider cooperation, since in the first place, there is no other conscious
agent to cooperate with. \\

And yet, in the typical practical situations when the single decision maker is facing multiple
and conflicting objectives, all the mentioned seeming advantages are instantly cancelled.
Instead, the single decision maker can easily end up feeling as if two or more autonomous
agents have moved inside of him or her, and now he or she has to turn into a dictator who, in
fact, ends up fighting himself or herself. \\
The obvious way out, therefore, for the single decision maker is to set up certain forms of
cooperation or bargaining between his or her multiple and conflicting objectives. And as seen
in section 6, this again is far from being a trivial task. Yet, the advantage remains that
there is only one single conscious and rational agent involved, namely, the decision maker,
and therefore, there is only one person to convince about the advantages of specific ways of
cooperations. \\

And in fact, in this case there is no need to convince anybody about the critical importance
and usefulness of cooperation. \\

In this way, cooperation can be seen as being highly pertinent to all the three theories of
games, choice and decision. And in some ways, cooperation may be more at home, even if not at
ease as well, in decisions. \\

What is suggested in this study is, therefore, to face consciously the following
alternative : \\

{\bf Non-cooperation} ~:~ Either we limit rationality and try to use it as much as possible
within the framework of non-cooperation, and only on occasion, and only as a second choice do
we try cooperation as well. \\

Or \\

{\bf Cooperation} ~:~ Within an extended and deeper sense of rationality, in fact, a
meta-rationality, we create and maintain a context in which the priority given to cooperation
can be relied upon. In other words, we are {\it firmly and reliably cooperation minded}. \\

So far, the first alternative has mostly been taken, even if rather by default. In games, as
follows from Binmore, this limitation of rationality mostly to the study of non-cooperative
games is still leaving us with complexities which are not algorithmically solvable. Thus,
having mostly avoided what appeared a theoretically difficult task in the second alternative,
instead of it, the first alternative, which in fact is not less difficult, was taken in game
theory. In social, collective or group choice, as seen in section 5, there are immense
possibilities for cooperation, while the alternative is the introduction of a dictator. And as
far as decision theory is concerned, what may be somewhat surprising, it proves to be a rather
natural and also unavoidable, even if not at all easy, framework for cooperation. \\

The way, therefore, suggested in this study is to take the second above alternative, namely,
of cooperation. \\
Here we have to keep in mind that both the non-cooperative and cooperative types of games are
of an extreme complexity from the point of view of theoretical approaches. And similar extreme
difficulties we face in social, collective or group choice, as well as in decisions. \\
Therefore, the focus can now shift away from the earlier attempts to construct comprehensive
enough theories of cooperation  for either of the three theories of games, choice or decisions.
Instead, one can focus on finding corresponding {\it cooperation minded frameworks} which are
large enough to contain much of what is already known and it is important about cooperation,
and then, in such frameworks, to develop further theories, methods, examples, applications,
and so on, related to cooperation, which have their own value and interest, even if they may
fall short of being comprehensive enough. \\

And clearly, once such a shift of focus is found appropriate, the {\it second alternative} is
the {\it natural} one. Indeed, this second alternative does obviously contain the first one as
a particular case. On the other hand, the first one  - in spite of the Nash Program, among
others - has so far not been proven to contain the second one. \\

So far, because of various reasons, not all of them fully conscious and rational, and some of
them related to the significant complexities involved in both theoretical and practical
aspects of cooperation, there has been a tendency to pursue the first above alternative.
However, as seen in this study, there are very strong reasons to opt firmly and reliably for
the second alternative above, namely, of cooperation. Indeed, in games, even the paradigmatic
non-cooperative Nash equilibrium loses all meaning outside of certain very strong cooperative
type assumptions. In social, collective or group choice, the failure to cooperate must bring
in a dictator. And in decision making with multiple conflicting objectives, the need for, and
advantages of cooperation do not have to be argued, the only issue being a satisfactory {\it
competence} in this regard. \\ \\

{\bf Note :} References [3], [4], [9], [11], [12], [15], [16], [22], [23], [28], [29], [33]
and [34] are particularly instructive for those who wish to become familiar with the subjects
in this study. References [11], [12], [15], [22] and [23] are still some of the very best of
their kind.


\begin{thebibliography}{99}

\vspace{0.5cm}

\bibitem{} Arrow, Kenneth J [1] : A difficulty in the concept of social welfare. Journal of
Political Economy, 58, 4 (1950)

\bibitem{} Arrow, Kenneth J [2] : Social Choice and Individual Values, 2nd ed. Wiley, New
York, 1963

\bibitem{} Axelrod, R : The Evolution of Cooperation. New York, Basic Books, 1984

\bibitem{} Bacharach, Michael \& Hurley Susan (Eds) : Foundations of Decision Theory,
Blackwell, Cambridge, 1991

\bibitem{} Binmore, Kenneth [1] : Modelling rational players. Part I. Economics and
Philosophy, 3 (1987) 179-214

\bibitem{} Binmore, Kenneth [2] : Modelling rational players. Part II. Economics and
Philosophy, 4 (1988) 9-55

\bibitem{} Binmore, Kenneth [3] : Game theory and the social contract : Mark II
(manuscript 1988) London School of Economics

\bibitem{} Blau, Julian : The existence of social choice functions. Econometrica, 25, 2
(1957) 302-313

\bibitem{} Hargreaves, Shaun P, et. al. : Game Theory, A Critical Intorduction. Routledge,
London, 1995

\bibitem{} Harsanyi, J C : Approaches to the bargaining problem before and after the theory
of games : a critical discussion of Zeuthen's, Hick's, and Nash's theories. Econometrica,
24 (1956), 144-157

\bibitem{} Kelly, Jerry S : Arrow Impossibility Theorems. Acad. Press, New York, 1978

\bibitem{} Luce, R Duncan \& Raiffa, Howard : Games and Decisions,
Introduction and Critical Survey. Wiley, New York, 1957, or Dover, New York, 1989

\bibitem{} Manin, Yu I : Mathematics and Physics. Birkhauser, Boston, 1981

\bibitem{} McKinsey, J C C : Introduction to the Theory of Games. Mc-Graw-Hill, New York, 1952

\bibitem{} Mirkin, Boris G : Group Choice. Wiley, New York, 1979

\bibitem{} Nasar, Silvia : A Beautiful Mind. Faber and Faber, London, 1998

\bibitem{} Nash, John F [1] : Equilibrium points in n-person games. Proc. Nat. Acad. Sci.
USA, 38 (1950), 48-49

\bibitem{} Nash, John F [2] : The bargaining problem. Econometrica, 18 (1950) 155-162

\bibitem{} Nash, John F [3] : Non-cooperative games. Ann. Math., 54 (1951) 286-295

\bibitem{} Nash, John F [4] : Two-person cooperative games. Econometrica, 21 (1953) 128-140

\bibitem{} von Neumann, John : Zur Theorie der Gesellschaftsspiele. Math. Annalen, 100
(1928) 295-320

\bibitem{} von Neuman, John \& Morgenstern, Oskar : Theory of Games and Economic
Behavior. Princeton, 1944

\bibitem{} Owen, Guillermo : Game Theory. Saunders, Philadelphia, 1968

\bibitem{} Pareto, Vilfredo : Course d'\'economie politique profess\'e \`a l'universit\'e de
Lausanne, 3 volumes, 1896-7

\bibitem{} Rasmusen, Eric : Games and Information. Balckwell, Malden, 2001

\bibitem{} Rosinger, Elemer E [1] : Interactive algorithm for multiobjective optimization.
JOTA, 35, 3 (1981) 339-365

\bibitem{} Rosinger, Elemer E [2] : Errata Corrige : Interactive algorithm for
multiobjective optimization. JOTA, 38, 1 (1982) 147-148

\bibitem{} Rosinger, Elemer E [3] : Aids for decision making with conflicting objectives.
In Serafini, P (Ed.), Mathematics of Multiobjective Optimization. Springer, New York,
1985, 275-315

\bibitem{} Rosinger, Elemer E [4] : Beyond preference information based multiple criteria
decision making. European Journal of Operational Research, 53 (1991) 217-227

\bibitem{} Tucker, Albert W : A two person dilemma. (unpublished) Stanford University
mimeos, May 1950.

\bibitem{} Vincke, Philippe [1] : Aggregation of preferences : a review. European Journal of
Operational Research, 9 (1982) 17-22

\bibitem{} Vincke, Philippe [2] : Arrow's theorem is not a surprising result. European
Journal of Operational Research, 10 (1982) 22-25

\bibitem{} Vorob'ev, N N : Game Theory, Lectures for Economists and Systems Scientists.
Springer, New York, 1978

\bibitem{} Walker, Paul : An outline of the history of game theory.
(http://william-king.www.drexel.edu/top/class/histf.html)

\bibitem{} Zeuthen, F : Problems of Monopoly and Economic Warfare. Routledge,
London, 1930 \\ \\

\end{thebibliography}
\end{document}